\documentclass[12pt]{article}

\usepackage{url}
\makeatletter
\newcommand{\address}[1]{\gdef\@address{#1}}
\newcommand{\email}[1]{\gdef\@email{\url{#1}}}
\newcommand{\sites}[1]{\gdef\@sites{\url{#1}}}
\newcommand{\@endstuff}{\par\vspace{\baselineskip}\noindent\small
	\begin{tabular}{@{}l}\scshape\@address\\\textit{E-mail address:} \@email \\ \textit{URL:} \@sites \end{tabular}}
\AtEndDocument{\@endstuff}
\makeatother

\title{Rational cross-sections, bounded generation \\ and orders on groups}
\author{Corentin Bodart}
\address{Section de Mathématiques, Université de Genève, Switzerland}
\email{corentin.bodart@unige.ch}
\sites{https://sites.google.com/view/corentin-bodart}

\usepackage{import}
\usepackage[backend=bibtex, style=alphabetic, sorting=nyt,
    doi=false,isbn=false,url=false,eprint=false, maxnames=50]{biblatex}
\usepackage[english]{babel}
\usepackage{amsmath}
\usepackage{amsfonts}
\usepackage{amssymb}
\usepackage{mathrsfs} % Cursive
\usepackage{dsfont}

\usepackage[dvipsnames]{xcolor}

\usepackage{hyperref}
\usepackage{enumitem}
\usepackage[chapter]{easy-todo}         % todo

\usepackage[left=3cm,right=3cm,top=2.8cm,bottom=3cm, marginparsep=3mm]{geometry}
\usepackage{changepage}

\usepackage[retainorgcmds]{IEEEtrantools}   % IEEEeqnarray
\usepackage[font=small, hypcap=false]{caption}
\usepackage{marginnote}
\usepackage{parskip}
\usepackage{tikz-cd}
\usepackage{tikz}
\usetikzlibrary{shapes, calc, arrows, patterns, decorations.markings, decorations.pathreplacing}
\usepackage{xifthen}
\pgfdeclarelayer{background} 
\pgfdeclarelayer{foreground}
\pgfsetlayers{background,main,foreground}
\import{Packages/}{custom_macros.tex}
\import{Packages/}{custom_theorem.tex}

\begin{document}

\maketitle

\begin{abstract}
    We provide new examples of groups without rational cross-sections (also called regular normal forms), using connections with bounded generation and rational orders on groups. Our examples contain a finitely presented HNN extension of the first Grigorchuk group. This last group is the first example of finitely presented group with solvable word problem and without rational cross-sections. It is also not autostackable, and has no left-regular complete rewriting system.
    \vspace*{2mm}
    
    \let\thefootnote\relax\footnotetext{The author acknowledges support of the Swiss NSF grants 200020-178828 and 200020-200400.} 
\end{abstract}
A rational cross-section for a group $G$ is a regular language $\Lc$ of unique representatives for elements of $G$. This notion can be traced back to Eilenberg and Schützenberger \cite{EILENBERG1969}, and was mainly explored by Gilman \cite{Gilman1987GROUPSWA}. A related \say{Markov property} was introduced by Gromov \cite{Gromov1987} and explored by Ghys and de la Harpe \cite{de1990groupes}.

Rational cross-sections are linked to several subjects in group theory. For instance constructing a geodesic rational cross-section (i.e., rational cross-section with minimal-length representatives) over a given generating set is sufficient to prove the rationality of the corresponding growth series. In particular, the following groups admit geodesic rational cross-sections for all generating sets:
\begin{itemize}[leftmargin=7mm, label=\textbullet]
    \item Finitely generated abelian groups \cite{Neumann_Shapiro}
    \item Hyperbolic groups \cite{Cannon, Gromov1987, de1990groupes}
\end{itemize}
If you don't require the section to be geodesic (and we won't), the property of having a rational cross-section is independent of the choice of a generating set. The following groups have been shown to admit rational cross-sections:
\begin{itemize}[leftmargin=7mm, label=\textbullet]
    \item Automatic groups \cite[Theorem 2.5.1]{Eps92}
    \item Finitely generated Coxeter groups \cite{Brink1993}
    \item Thompson's group $F$ \cite{guba} \vspace*{1mm}
\end{itemize}

We should also mention that the class of groups with rational cross-sections is closed under commensurability, extensions, free products \cite{Gilman1987GROUPSWA}, graph products \cite{Hermiller_GraphProduct} and some Bass-Serre constructions \cite{Hermiller2018}. As a corollary,
\begin{itemize}[leftmargin=7mm, label=\textbullet]
	\item virtually polycyclic groups (including virtually nilpotent groups),
	\item right-angled Artin groups (RAAG's) and
	\item Baumslag-Solitar groups $BS(m,n)$
\end{itemize}
all admit rational cross-sections.
\medbreak
On the other side, very few groups are known not to have rational cross-sections:
\begin{itemize}[leftmargin=7mm, label=\textbullet]
    \item Infinite torsion groups (see \cite{Gilman1987GROUPSWA, de1990groupes})
    \item Recursively presented groups with undecidable word problem (see eg.\ \cite{OTTO1998621}).
\end{itemize}
For completeness, we can also construct new non-examples from old ones:
\begin{thm*}[{\cite{GerstenShort}, see also Theorem \ref{sec2:nightmare}}]
Consider a product $A*_CB$ with $C$ finite, s.t.\ either $A$ or $B$ doesn't have a rational cross-section. Then neither does $A*_CB$.
\end{thm*}
\medbreak
Our main objective in the present paper is to add to the list of known non-examples. We start with the observation (Corollary \ref{sec3:dichotomy}) that groups with rational cross-sections either have bounded generation (see \S3 for a definition) or contain free non-abelian submonoids. Some new non-examples follow, specifically
\begin{itemize}[leftmargin=7mm, label=\textbullet]
    \item some extensions of infinite torsion groups (\S 3.1)
    \item some Grigorchuk-type groups (including the Fabrykowski-Gupta group) (\S 3.2)
\end{itemize}
We then turn our attention to wreath products. We define property (R+LO), a strengthening of left-orderability (LO), which means that some positive cone admits a rational cross-section. We derive the following (positive) result:
\begin{thm}[Theorem \ref{sec4:positive_result}]
Let $L$ and $Q$ be two groups such that
\begin{itemize}[leftmargin=7mm, label=\upshape\textbullet]
    \item $L$ has a rational cross-section
    \item $Q$ is virtually {\upshape(R+LO)}
\end{itemize}
then the restricted wreath $L\wr Q$ has a rational cross-section.
\end{thm}
In the other direction, we prove that $Q$ admitting \say{large} rational partial orders is actually a necessary condition for $L\wr Q$ to admit a rational cross-section (under extra conditions on $L$). Such rational orders have received attention recently (see \cite{hermiller2016positive, su2020formal, alonso2020geometry, Antoln2021RegularLO}). We build on these results to prove

\newpage

\begin{thm}
Let $Q$ be a group satisfying either
\begin{enumerate}[leftmargin=8mm, label={\upshape(\alph*)}]
    \item $Q$ has arbitrarily large/infinite torsion subgroups; \hfill {\upshape(Theorem \ref{sec5:large_torsion})}
    \item $Q$ has infinitely many ends; \hfill {\upshape(Theorem \ref{sec5:infinite_ended})}
\end{enumerate}
and let $L$ be a non-trivial group without free non-abelian submonoid, then $L\wr Q$ does not admit any rational cross-section.
\end{thm}

This result provides many \textit{elementary} non-examples. For instance, Ghys and de la Harpe ask whether all finitely generated solvable groups have rational cross-sections. Kharlampovich  already constructed finitely presented $3$-step solvable groups with unsolvable word problem \cite{3Step-Solvable}, which therefore cannot admit any rational cross-section. To the best of our knowledge, groups like $C_2\wr (C_2\wr \ZZ)$ covered in case (a) are the first set of examples of $3$-step solvable groups with solvable word problem and without rational cross-sections. Examples like $\ZZ\wr F_2$ covered in case (b) give the first torsion-free non-examples which \say{do not rely on the existence of infinite torsion groups}.

Finally, we prove a bounded-generation-like condition for torsion-by-$\ZZ$ groups admitting rational cross-sections (Proposition \ref{sec6:torsion-by-Z}). This result (and some more work) allows us to prove the following two theorems:

\begin{thm}[Theorem \ref{sec6:Houghton_is_down}]
The permutation group
\[ H_2 = \left\{ \sigma \in \Sym(\ZZ) \;\Big|  \begin{array}{c} \exists \pi\in \ZZ \text{ such that }\sigma(x)=x+\pi \\ \text{ for all but finitely many }x \in \ZZ \end{array} \right\} \]
doesn't have any rational cross-section.
\end{thm}
This group is often called \textit{second Houghton group} as part of the family introduced in \cite{Houghton1978TheFC}. Note that $H_2=\la a,t\ra$ with $a=(1\,2)$ and $t(x)=x+1$. We should mention B.\ Neumann introduced an alternating version $\la b,t\ra$ with $b=(1\,2\,3)$ much earlier \cite{Neumann}.
\begin{thm}[Theorem \ref{sec7:Lysenok_is_down}]
The following HNN extension of the first Grigorchuk group
\[ \Lk = \la a,b,c,d,t \;\;\Big| \begin{array}{c} a^2=b^2=c^2=d^2=bcd=(ad)^4=(adacac)^4 = e \\ t^{-1}at = aca,\; t^{-1}bt=d,\;t^{-1}ct=b,\; t^{-1}dt=c \end{array}\ra \]
(introduced in \cite{Grigorchuk_1998}) doesn't have any rational cross-section.
\end{thm}
This last example is the first known example of finitely presented group with solvable word problem and without rational cross-section, answering \cite[Q4]{OTTO1998621}. As a corollary, it also provides the first example of finitely presented group with solvable word problem that is not autostackable, answering a question of Hermiller. Another corollary is that it is not presented by a left-regular complete rewriting system.
\medbreak
\textbf{Acknowledgment}.\ I would like to thanks Tatiana Nagnibeda for introducing me to this problem, and for \textit{many} fruitful discussions. I'm also grateful to Susan Hermiller and Laurent Bartholdi for sharing some of their thoughts and ideas on the issue. I'd like to thank Rostislav Grigorchuk and Dmytro Savchuk for pointing to a mistake in the proof of Proposition \ref{sec7:action_on_ternary} in a previous version of this paper.

\newpage

\counterwithin{thm}{section}

\section{Background information}
\subsection{Regular languages}
Regular languages form the lowest class of languages in the Chomsky hierarchy of complexities. Informally, a language is \textit{regular} if its membership problem can be decided by some computer with finite memory. % \footnote{This is really the definition of "recognizable subset of $\Ac^*$", but "rational" and "recognizable" are equivalent in the free monoid $\Ac^*$}
Our model of computer is the following:
\begin{defi}
An \textbf{automaton} is a $5$-uple $M=(V,\Ac,\delta,*,T)$ where
\begin{itemize}[leftmargin=8mm]
    \item $V$ is a set of states / vertices.
    \item $\Ac$ is an alphabet.
    \item $\delta\subseteq V\times\Ac\times V$ is the transition function. An element $(v_1,a,v_2)\in \delta$ should be seen as an oriented edge from $v_1$ to $v_2$, labeled by $a$.
    \item $*\in V$ is the initial vertex.
    \item $T\subseteq V$ is the set of \say{accept} / terminal vertices.
\end{itemize}
An automaton is \textbf{finite} if both $V$ and $\delta$ are finite.
\end{defi}

Here are some examples of finite automata (with terminal vertices in green):

\begin{center}
\begin{minipage}{.6\linewidth}
\centering
\vspace*{2mm}
	\begin{tikzpicture}[scale=1.1, thick]
		\node[circle, draw=Green, ultra thick, inner sep=1pt] (v0) at (0,0) {$*$};
		\node[circle, fill=Green, inner sep=2pt, outer sep= 3pt] (v+) at (2,0) {};
		\node[circle, fill=Green, inner sep=2pt, outer sep= 3pt] (v-) at (-2,0) {};
		
		\begin{scope}[every node/.style={fill=white, inner sep=1.7pt}]
		    \draw[-latex] (v0) to node[pos=.45, inner sep=2.5pt]{$t$} (v+);
		    \draw[-latex] (v0) to node[pos=.45, inner sep=2.5pt]{$T$} (v-);
		    
		    \draw[-latex, out=30, in=-30, looseness=16] (v+) to node[pos=.35]{\small$t$} (v+);
		    \draw[-latex, out=150, in=210, looseness=16] (v-) to node[pos=.35]{\small$T$} (v-);

		\end{scope}
	\end{tikzpicture}
    \captionof{figure}{$M$ for $\Lc=\{\varepsilon\}\cup\{t^n, T^n\mid n\ge 1\}$ \label{auto:Z}}
\vspace*{8mm}
    \begin{tikzpicture}[scale=1.25, thick]
		\node[circle, draw=black, ultra thick, inner sep=1pt] (v0) at (-1,0) {$*$};
		\node[circle, fill=Green, inner sep=2pt, outer sep= 3pt] (v1) at (.5,.87) {};
		\node[circle, fill=Green, inner sep=2pt, outer sep= 3pt] (v2) at (.5,-.87) {};
		\node[circle, fill=black, inner sep=2pt, outer sep= 3pt] (v3) at (2,-.5) {};
		\node[circle, fill=black, inner sep=2pt, outer sep= 3pt] (v4) at (2,.5) {};
		
		\begin{scope}[every node/.style={fill=white, inner sep=1.7pt}]
		    \draw[-latex] (v0) to node[pos=.45, inner sep=2.5pt]{$t$} (v1);
		    \draw[-latex] (v1) to node[pos=.45, inner sep=2.5pt]{$t$} (v2);
		    \draw[-latex] (v2) to node[pos=.45, inner sep=2.5pt]{$t$} (v0);
		    \draw[-latex] (v2) to node[pos=.45, inner sep=2.5pt]{$T$} (v3);
		    \draw[-latex] (v4) to node[pos=.45, inner sep=2.5pt]{$t$} (v1);
		\end{scope}
	\end{tikzpicture}
	\captionof{figure}{$M$ for $\Lc=\{t^n\mid n\equiv 1,2\pmod 3\}$ \label{auto:mod3}}
\end{minipage}
\begin{minipage}{.38\linewidth}
\centering
    \begin{tikzpicture}[thick, rotate=-90]
\node[circle, draw=Green, ultra thick, inner sep=1pt] (0) at (0,0) {$*$};
\node[circle, fill=Green, inner sep=2pt, outer sep= 3pt] (p1) at (1,1.5) {};
\node[circle, fill=Green, inner sep=2pt, outer sep= 3pt] (m1) at (1,-1.5) {};
\node[circle, fill=Green, inner sep=2pt, outer sep= 3pt] (a1) at (2.3,0) {};
\node[circle, fill=black, inner sep=2pt, outer sep= 3pt] (p2) at (4.7,0) {};
\node[circle, fill=Green, inner sep=2pt, outer sep= 3pt] (die+) at (3.3,1.5) {};
\node[circle, fill=Green, inner sep=2pt, outer sep= 3pt] (die-) at (3.3,-1.5) {};

\begin{scope}[every node/.style={fill=white, inner sep=1.7pt}]
    \draw[-latex] (0) -- (p1) node [pos=.45] {$t$};
    \draw[-latex] (0) -- (m1) node [pos=.45] {$T$};

    \draw[-latex, red] (0) -- (a1) node[pos=.4] {$a$};
    \draw[-latex, red] (p1) -- (a1) node[pos=.4] {$a$};
    \draw[-latex, red] (m1) -- (a1) node[pos=.4] {$a$};

    \draw[-latex, bend left=20] (a1) to node[pos=.5]{$t$} (p2);
    \draw[-latex, red, bend left=20] (p2) to node[pos=.5]{$a$} (a1);
    
    \draw[-latex] (a1) to node[pos=.5]{$T$} (die-);
    \draw[-latex] (a1) to node[pos=.45]{$t$} (die+);
    
    \draw[-latex, out=50, in=110, looseness=12] (p1) to node[pos=.5]{\small$t$} (p1);
    \draw[-latex, out=-70, in=-130, looseness=12] (m1) to node[pos=.5]{\small$T$} (m1);
    
    \draw[-latex, out=22, in=-38, looseness=12] (p2) to node[pos=.5]{\small$t$} (p2);
    \draw[-latex, out=50, in=110, looseness=12] (die+) to node[pos=.5]{\small$t$} (die+);
    \draw[-latex, out=-70, in=-130, looseness=12] (die-) to node[pos=.5]{\small$T$} (die-);
\end{scope}
\end{tikzpicture}
    \captionof{figure}{$M$ for $\Lc\underset{\ev}{\longleftrightarrow} C_2\wr \ZZ$ \label{auto:wr}}
\end{minipage}
\end{center}

\begin{defi}
A finite automaton $M$ \textbf{recognizes} a word $w\in\Ac^*$ if $w$ can be read along some oriented path from the initial state $*$ to some terminal state $v\in T$. A language $\Lc\subseteq \Ac^*$ is \textbf{regular} if it can be written as
\[ \Lc = \{ w\in \Ac^* \mid w\textnormal{ is recognized  by }M \} \]
for some finite automaton $M=(V,\Ac,\delta,*,T)$.
\end{defi}
\textbf{Remark.} Whenever $*\in T$ the automaton $M$ recognizes the empty word (noted $\varepsilon$).

\bigbreak

Some additional terminology around automata will be needed:
\begin{defi} $\,$
\begin{itemize}[leftmargin=8mm]
    \item An automaton is \textbf{deterministic} whenever, for all $v\in V$ and $a\in \Ac$, there exists at most one edge exiting $v$ and labeled by $a$ (i.e., $\abs{\big(\{v\}\times\{a\}\times V\big)\cap \delta}\le 1$).
    
    \item For $v_1,v_2\in V$, we say that $v_2$ is \textbf{accessible} from $v_1$ is there exists an \textit{oriented} path from $v_1$ to $v_2$ in $M$. An automaton is  \hypertarget{trimmed}{\textbf{trimmed}} if, for all $v\in V$, there exists $t\in T$ such that $t$ is accessible from $v$ and $v$ is accessible from $*$.
    
    \item A \textbf{strongly connected component} is a maximal subset $C\subseteq V$ such that, for every $v_1,v_2\in C$, the state $v_2$ is accessible from $v_1$ (and reciprocally).
\end{itemize}
\end{defi}
For instance, automata in Figure \ref{auto:Z} and \ref{auto:mod3} are deterministic, but not in Figure \ref{auto:wr}. The automaton in Figure \ref{auto:mod3} is not trimmed.

A remarkable result due to Rabin and Scott \cite{RabinScott} is that any regular language is the set of words recognized by some finite, \emph{deterministic}, trimmed automaton. Another fundamental result in the theory of regular language is Kleene's theorem. (We will use it without even mentioning it, usually to construct regular languages efficiently.)

\begin{thm}[Kleene's theorem]
The class of regular languages (over a finite alphabet $\Ac$) is the smallest class of languages containing finite languages (over $\Ac$), and closed under the following three operations:
\begin{itemize}[leftmargin=8mm]
    \item Finite union
    \item Concatenation: $\Lc_1\Lc_2=\{w_1w_2\mid w_1\in\Lc_1,w_2\in\Lc_2\}$.
    \item Kleene's star: $\Lc^*=\bigcup_{n\ge 0}\Lc^n=\{w_1w_2\ldots w_n\mid w_i\in\Lc, n\ge 0\}$
\end{itemize}
The class of regular languages is also closed under complementation (in $\Ac^*$), hence finite intersection, set difference, and in general all Boolean operations.
\end{thm}

\medbreak

Finally, we introduce some notations for subwords (in general languages).
\begin{nota} \label{sec1:nota}
Consider a word $w\in\Ac^*$. We denote by $\ell=\ell(w)$ its length. For $0\le i< j\le \ell$, we denote by $w(i\col j]$ the subword consisting of all letters from the $(i+1)$-th to the $j$-th one (included). If $i=0$ we abbreviate the notation to $w(\icol j]$.
\end{nota}

%Finally, we'll add some terminology on general languages :
%\begin{defi}$\,$
%\begin{itemize}[leftmargin=8mm]
%    \item A language $\Lc$ is \textbf{prefix-closed} if $uv\in\Lc$ implies $u\in \Lc$. 
%    \item A language $\Lc$ is \textbf{$+$-prefix-closed} if $uv\in\Lc$ and $u\ne\varepsilon$ implies $u\in\Lc$.
%\end{itemize}
%A regular language is prefix-closed if and only if it is recognized by some finite automaton with only  terminal vertices (i.e. $T=V$).
%\end{defi}
%\bigbreak

%Languages associated to automata in Figure \ref{auto:Z} and \ref{auto:wr} are prefix-closed.

\subsection{Rational subsets and cross-sections}
We now export those definitions to groups
\begin{defi}
Let $G$ be a group.
\begin{itemize}[leftmargin=8mm]

    \item Given a generating set $\Ac$ and a word $w\in\Ac^*$, this word can be evaluated in $G$. The corresponding element of $G$ will be denoted $\ev(w)$ or $\overline w$.
    
    \item A subset $R\subseteq G$ is \textbf{rational} if there exists a generating set $\Ac$ and a regular language $\Lc\subseteq\Ac^*$ s.t.\ $\ev(\Lc)=R$. We denote the class of rational subsets of $G$ by $\Rat(G)$. For example finitely generated subgroups $H\le G$ are rational.
    
    \item If moreover the evaluation map $\ev\colon\Lc\to R$ is bijective, we say that $\Lc$ is a \textbf{rational cross-section} for $R$.
    %\item If moreover $\Lc$ is prefix-closed, we say $\Lc$ is a \textbf{Markov language} for $R$.
\end{itemize}
\end{defi}
Our primary interest will be rational cross-sections for the entire $R=G$.

\medbreak

For instance, automaton in Figure \ref{auto:Z} recognizes a rational cross-section for $\ZZ=\la t\ra$ relative to the generating set $\Ac=\{t,T=t^{-1}\}$, while automaton in Figure \ref{auto:wr} recognizes a rational cross-section for $C_2\wr \ZZ=\la a\ra\wr\la t\ra$ (with $\Ac=\{a,t,T=t^{-1}\}$).

\bigbreak

It is interesting to know which operations preserve rationality of subsets of a given group $G$. Several results directly translate from regular languages to rational subsets: the class of rational subsets of a group $G$ is closed under finite union, set-theoretic product, Kleene star. However, contrary to intersections of regular languages, intersections of rational subsets are not rational in general. That being said, things behave better in special cases. Let us for instance recall the following classical result (used in \S 5).
\begin{lemma}[Benois' lemma {\cite{Benois1}, see also \cite{Rational_Subset}}] \label{Benois}
Consider $F_2$ the free group of rank $2$, and let $R\subseteq F_2$ be a rational subset. Then the set of reduced words representing elements of $R$ is regular. As a corollary, $\Rat(F_2)$ is closed under all Boolean operations.
\end{lemma}

\medbreak

A natural question is how  properties \say{rationality} and \say{existence of rational cross-section} for $R$ depend on the choice of a monoid generating set $\Ac$ and on the ambient group $G$. It is answered by Gilman \cite{Gilman1987GROUPSWA}:

\begin{prop} \label{sec1:dependance} Let $G$ be a group and $R\subseteq G$ a rational subset,
\begin{enumerate}[leftmargin=8mm, label={\upshape(\alph*)}]
    \item Let $\Ac$ be any set generating $G$ as a monoid. There exists regular language $\Lc\subseteq\Ac^*$ such that $\ev(\Lc)= R$. If furthermore $R$ has a rational cross-section, we can ensure $\ev\colon\Lc\to R$ is bijective (possibly for another $\Lc\subseteq\Ac^*$).
    \item The subgroup $\la R\ra$ is finitely generated, and $R$ is rational in $\la R\ra$. It follows that $R$ is rational in any ambient group $H\supseteq R$. Moreover, the same holds with \say{rational} replaced by \say{admitting a rational cross-section} everywhere. 
\end{enumerate}
\end{prop}

%Note that the dependence of "existence of a Markov language for $R=G$" on the choice of a generating set remains open. However it is known that, if $G$ has a rational cross-section, then it has a Markov language for \textbf{some} generating set $\Ac$ (see \cite{Gilman1987GROUPSWA}).
\subsection{Left-invariant orders}

Let $G$ be a group. A \textbf{left-invariant order} on $G$ is a (partial) order $\prec$ on $G$ satisfying 
\[ h_1\prec h_2\implies gh_1\prec gh_2 \quad\text{for all }g,h_1,h_2\in G. \]
Given a left-invariant order on $G$, we can define its \textbf{positive cone}
\[ P_\prec = \{ g\in G \mid g\succ e_G \}. \]
Note that $P_\prec^{-1}=\{g\in G\mid g\prec e_G\}$. Each property of the relation $\prec$ translates into a property of its positive cone,
\begin{enumerate}[leftmargin=7mm, label=(\alph*)]
    \item Anti-symmetry translates into $P_\prec \cap P_\prec^{-1}=\emptyset$.
    \item Transitivity translates into the fact that $P_\prec$ is a sub-semigroup, i.e., $P_\prec P_\prec\subseteq P_\prec$.
    \item The order $\prec$ is total if and only if $G=P_\prec\cup \{e\}\cup P_\prec^{-1}$.
\end{enumerate}
The other way around, given a subset $P\subseteq G$ satisfying both (a) and (b), we can define a left-invariant (partial) order
\[ g\prec_P h\iff g^{-1}h\in P.\]
Hence we can think interchangeably about orders and their positive cones. This allows us to define a \textbf{rational order} (sometimes \textbf{regular order}) on $G$ as a left-invariant order on $G$ whose positive cone is a rational subset of $G$.

\textbf{Remark.} Most of the recent literature deals with \textbf{total} left-invariant orders, so that this adjective is usually dropped. We will work with both, hence keep the adjective.

\section{Intersections of rational subsets}

In the spirit of Benois' lemma, we study when rationality is preserved when taking intersections of subsets. First, we have

\begin{prop}[{Compare with \cite[Proposition 4.5]{GerstenShort}}]\label{sec2:nightmare}
Let $G$ be a group, $H$ a subgroup, and $R\subseteq G$ a rational subset. If either
\begin{enumerate}[leftmargin=8mm, label={\upshape(\alph*)}]
    \item $G=H*_CB$ over a finite subgroup $C$, or
    \item $G=H*_Ct$ over finite subgroups $\iota,\iota':C\into H$,
\end{enumerate}
then $H\cap R$ is rational. Moreover, if $R$ had a rational cross-section, so will $H\cap R$.
\end{prop}
\textbf{Remark.} By induction, the same result holds for $G$ the fundamental group of a graph of groups with finite edge groups, and $H$ any vertex group.

\begin{proof} We will only prove part (a) as part (b) is similar.

Let us first treat two small cases: If $C=B$ then $H=G$ hence $H\cap R=R$ is still rational. If $C=H$, then $H$ is finite hence $H\cap R$ is rational.

Otherwise we can suppose $\Ac\subseteq (H\cup B)\setminus C$. Let $\Lc\subseteq \Ac^*$ be a regular language for $R$, and $M=(V,\Ac,\delta,*,T)$ an automaton recognizing $\Lc$. We construct a new automaton  $M'=(V,(\Ac\setminus B)\sqcup C,\delta',*,T)$ as follows:
    \begin{itemize}[leftmargin=5mm]
        \item For each pair of states $p,q\in V$, if there exists an (oriented) path from $p$ to $q$ starting with an edge labeled by $s\in B$, with associated word evaluating to $c\in C$ and no proper prefix evaluating in $C$, then add a $c$-edge from $p$ to $q$.
        
        \item Remove all edges labeled by letters $s\in B\setminus C$.
    \end{itemize}
    It should be clear that the language $\Lc'$ recognized by $M'$ evaluate to $H\cap R$. Let us now suppose $\Lc$ was a rational cross-section and prove that the same holds for $\Lc'$.
    
    Each word $w\in \Lc'$ can be written as $w=u_0c_1u_1\ldots c_m u_m$ with $u_i\in(\Ac\setminus B)^*$ and $c_i\in C$. Those words corresponds to words $w_\old\in \Lc$ of the form
    \[ w_\old = u_0v_1u_1\ldots v_mu_m \]
    with $v_i\in\Ac^*$ evaluating to $c_i$ (in particular $\bar w=\bar w_\old$), starting with a letter in $B$ and without proper prefix evaluating in $C$. Suppose we have two words $w,\tilde w\in\Lc'$ evaluating to the same element $g\in H\cap R$. It follows that both words $w_\old, \tilde w_\old$ evaluate to $g$, hence they are actually the same word:
    \[ u_0v_1\ldots v_mu_m = w_\old = \tilde w_\old = \tilde u_0\tilde v_1\ldots \tilde v_n\tilde u_n \]
    Comparing both expressions we get $u_0=\tilde u_0$ as those are the longest prefixes without letters in $B$. As we keep reading, the shortest word between $v_1$ and $\tilde v_1$ is a prefix of the other, but they both evaluate in $C$, hence $v_1=\tilde v_1$ and $c_1=\tilde c_1$. Iterating, we get $w=\tilde w$ as wanted: $\Lc'$ is a rational cross-section.
\end{proof}
Surprisingly, the reciprocal \say{if $A$ and $B$ have rational cross-sections (and $C$ is finite) then $A*_CB$ has a rational cross-section} seems open. Note that, using Proposition 3.3 from \cite{Hermiller2018}, this reduces to the following question:

\begin{adjustwidth}{3mm}{3mm}
    \textbf{Question.} Suppose that $A$ has a rational cross-section, and $C\le A$ is a finite subgroup. Is it true that $C$ has a regular language of cosets representatives? % Geometrically, if a (rooted, deterministic) labeled graph has a finite cover which has rational cross-section, does it imply the graph itself has a rational cross-section?
\end{adjustwidth}

Next we prove the following result. Note that Proposition \ref{sec2:nightmare} and Theorem \ref{sec2:nightmare_bis} lie at the opposite sides of a spectrum: the subgroup $H$ in Proposition \ref{sec2:nightmare} is a free factor, while here $H$ intersects each free factors (and each conjugate) at most once.

\begin{thm}\label{sec2:nightmare_bis} Let $G$ be a finitely generated group and $R$ a rational subset. Suppose
\begin{enumerate}[leftmargin=8mm, label={\upshape(\alph*)}]
    \item $G=A*_CB$ for $C$ a finite subgroup, or
    \item $G=A*_C t$ over finite subgroups $\iota,\iota':C\into A$,
\end{enumerate}
Suppose that $H\leqslant G$ is a finitely generated subgroup acting freely on the associated Bass-Serre tree (by left multiplication). Then $H\cap R$ is rational.
\end{thm}
\textbf{Remark.} These results can be made effective as soon as the rational subset membership problem is decidable in $G$ (i.e., is decidable in both factors, see \cite{KAMBITES2007622}).

We start with several lemmas, following the scheme in \S 3-4 of \cite{su2020formal}:
\begin{defi}
Let $G$ be a group, $\Ac$ a symmetric generating set, and $K\ge 0$ a constant. Consider $v,w\in\Ac^*$ two words. We say \textbf{$w$ asymmetrically $K$-fellow travel} with $v$ if there exists a weakly increasing sequence $(i_n)_{n=1,\ldots,\ell(w)}$ such that (using Notation \ref{sec1:nota})
\[ \dist_\Ac\big(\bar w(\icol n],\bar v(\icol i_n]\big) \le K \qquad\text{for all }n. \]

\end{defi}
\begin{center}
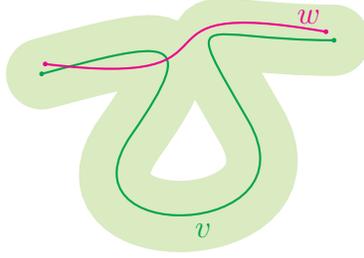

		\begin{tikzpicture}[scale=.7, rotate=-20, line cap=round]
			
			\clip[rotate=20] (0,-3.3) rectangle (7.7,2);
			\draw[draw=white, double=LimeGreen!30, thick, double distance=27pt, rounded corners=2]  (1,.5)
				to [out=35, in=75, looseness = 3]  (3,0) 
				to [out=-105, in=-40, looseness = 3] (5,1)
				to [out=140, in=-160, looseness = 3] (6,3);
				
			\draw[thick, Green]  (1,.5)
				to [out=35, in=75, looseness = 3]  (3,0) 
				to [out=-105, in=-40, looseness = 3] (5,1)
				to [out=140, in=-160, looseness = 3] (6,3);
				
			\draw[thick, Magenta]  (1,.7)
				to [out=13, in=190, looseness = 2]  (5.8,3.1);
				
			\draw[draw=Green, fill=Green] (1,.5) circle (1pt);
			\draw[draw=Green, fill=Green] (6,3) circle (1pt);
			
			\draw[draw=Magenta, fill=Magenta] (1,.7) circle (1pt);
			\draw[draw=Magenta, fill=Magenta] (5.8,3.1) circle (1pt);
			
			\node[Magenta] at (5.4,3.25) {$w$};
			\node[Green] at (4.9,-1.25) {$v$};
		\end{tikzpicture}
		\captionsetup{justification=centering,margin=1cm}
		\captionof{figure}{$w$ asymmetrically $K$-fellow travel with $v$. We see $w$ is strongly \\ restrained by $v$ while conditions imposed on $v$ by $w$ are much weaker.}
\end{center}

\begin{lemma}[{Compare with \cite[Lemma 3.5]{su2020formal}}] \label{sec2:su_lemma_35}
Let $G$ be a f.g. group, $\Ac$ a finite symmetric generating set, and $K\ge 0$. Consider $\Lc\subseteq\Ac^*$ a regular language. We define
\[ \tilde\Lc = \{ w\in\Ac^* \mid \exists v\in\Lc,\; \bar w=\bar v\text{ \upshape{and} }w\text{ \upshape{asymmetrically} }K\text{\upshape{-fellow travel with} }v\}. \]
Then $\tilde\Lc$ is rational and $\ev(\tilde\Lc)=\ev(\Lc)$.
\end{lemma}
\begin{proof}
Let $M=(V,\Ac,\delta,*,T)$ be a deterministic automaton recognizing $\Lc$. We construct a new automaton $\tilde M$ recognizing $\tilde \Lc$ as follows
\begin{itemize}[leftmargin=8mm, label=\textbullet]
    \item The vertex set is $\tilde V=V\times B_K$ where $B_K$ is the ball of radius $K$ (around $e$) in the Cayley graph of $G$, relative to the generating set $\Ac$.
    \item Add an $a$-edge ($a\in\Ac$) from $(p,g)$ to $(q,h)$ if any of the following conditions hold
    \begin{enumerate}[leftmargin=7mm, label=\arabic*.]
        \item $g=h=e$ and there was an $a$-edge from $p$ to $q$ in $M$,
        \item $p=q$ and $h=ga$ (i.e., there was an $a$-edge from $g$ to $h$ in $B_K$), or
        \item Using only edges of type 1 and 2, there already exists in $\tilde M$ an oriented path from $(p,g)$ to $(q,h)$ with associated word evaluating to $a$.
    \end{enumerate}
    \item The initial and terminal vertices are $\tilde *=(*,e)$ and $\tilde T=T\times\{e\}$ respectively.
\end{itemize}
Note that $M$ can be seen as the sub-automaton consisting of all type $1$ edges.
\begin{center}
\begin{tikzpicture}[scale=1,
    type1/.style={fill=black, inner sep=1.5pt, outer sep=2pt, circle},
    type2/.style={fill=Purple, inner sep=1pt, outer sep=1.5pt, circle},
    typeA/.style={circle, fill=white, inner sep=.8pt}]
    
    \node[type1] (1) at (0,0) {};
    \node[type1] (2) at (2,-1.5) {};
    \node[type1] (3) at (4,0) {};
    \node[type1] (4) at (6.5,1) {};
        
    \begin{scope}[very thick]
        \draw[-latex] (1) -- (2) node[pos=.5, circle, fill=white, inner sep=1pt] {$a$};
        \draw[-latex] (1) -- (3) node[pos=.5, circle, fill=white, inner sep=1pt] {$b$};
        \draw[-latex] (2) -- (3) node[pos=.5, circle, fill=white, inner sep=1pt] {$a$};
        \draw[-latex] (3) -- (4) node[pos=.5, circle, fill=white, inner sep=1pt] {$A$};
    \end{scope}
    
    \footnotesize{
    \begin{scope}[-latex, Purple]

            \node[type2] (11) at (0.3,0.9) {};
            \node[type2] (12) at (-0.3,-0.9) {};
            \node[type2] (13) at (-1,0.2) {};
        
        \begin{scope}[shift={(4,0)}]
            \node[type2] (31) at (-0.3,0.9) {};
            \node[type2] (32) at (0.3,-0.9) {};
            \node[type2] (33) at (1,-0.3) {};
        \end{scope}
        
        \begin{scope}[shift={(6.5,1)}]
            \node[type2] (41) at (-0.3,0.9) {};
            \node[type2] (42) at (0.3,-0.9) {};
            \node[type2] (43) at (1,-0.1) {};
        \end{scope}
        
        \begin{scope}[bend right]
            \draw (1) to node[pos=.47, typeA]{$a$} (11);
            \draw (11) to node[pos=.45, typeA]{$A$} (1);
            \draw (12) to node[pos=.47, typeA]{$a$} (1);
            \draw (1) to node[pos=.45, typeA]{$A$} (12);
        
            \draw (3) to node[pos=.47, typeA]{$a$} (31);
            \draw (31) to node[pos=.45, typeA]{$A$} (3);
            \draw (32) to node[pos=.47, typeA]{$a$} (3);
            \draw (3) to node[pos=.45, typeA]{$A$} (32);
        
            \draw (4) to node[pos=.47, typeA]{$a$} (41);
            \draw (41) to node[pos=.45, typeA]{$A$} (4);
            \draw (42) to node[pos=.47, typeA]{$a$} (4);
            \draw (4) to node[pos=.45, typeA]{$A$} (42);
        \end{scope}
        
        \draw[latex-latex] (1) to node[pos=.5, fill=white, inner sep=1.5pt]{$b$} (13);
        \draw[latex-latex] (3) to node[pos=.5, fill=white, inner sep=1.5pt]{$b$} (33);
        \draw[latex-latex] (4) to node[pos=.5, fill=white, inner sep=1.5pt]{$b$} (43);
        \end{scope}
        
        \draw[magenta, -latex, rounded corners=8pt] (1) .. controls (.8,-1) .. (2,-1.7) .. controls (3.2,-1) .. (4,0) node[pos=.25, typeA] {$a$} -- (32); 
        \draw[magenta, -latex, rounded corners=8pt] (11) -- (0,0) .. controls (2,.3) .. (4,0) node[pos=.6, fill=white, inner sep=1.5pt] {$b$} -- (31) ;
        \draw[magenta, -latex, rounded corners=8pt] (33) to[bend right=20] (4,0) .. controls (5.2,.85) .. (6.5,1) node[pos=.57, fill=white, inner sep=1.5pt] {$b$} -- (41) ;}
        
\end{tikzpicture}
\captionsetup{justification=centering,margin=1cm}

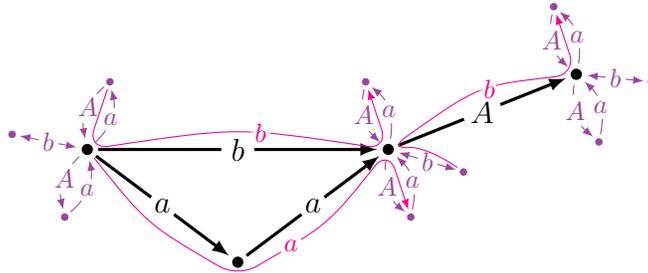
\captionof{figure}{Part of $\tilde M$ for $G=\la a,b\mid [a,b],b^2\ra$,  $\Ac=\{a,A=a^{-1},b\}$ and $K=1$. Type 1 edges can be seen in black, type 2 in purple, and some type 3 in pink.}
\end{center}

\newpage

We show that $\tilde\Lc$ coincide with the language recognized by $\tilde M$:
\begin{itemize}[leftmargin=8mm]
	\item By construction, each edge $e$ of type $3$ (from $(p,g)$ to $(q,h)$, labeled by $a\in\Ac$) can be associated to a word $v_e\in \Lc_{p\to q}$ satisfying $g^{-1}\bar v_e h=a$.

	Given a word $w$ recognized by $\tilde M$, we construct $v$ recognized by $M$ as follows: pick an accepted path for $w$ in $\tilde M$. Replace each edge of type $1$ by its label, forget about edges of type $2$, and replace each edge $e$ of type $3$ by $v_e$. It's easy to check that $\bar w=\bar v$ and $w$ asymmetrically $K$-fellow travel with $v$, hence $w\in\tilde\Lc$.

	\item Reciprocally, let $w\in\Ac^*$ be a word asymmetrically $K$-fellow traveling with some word $v\in \Lc$ evaluating to the same $\ev(w)=\ev(v)$. Pick an accepted path for $v$ in $M$. Each prefix $v(\icol i]$ leads to a state $p_i$ in $M$. Moreover, for each $0\le n\le \ell(w)$, there exists $g_n\in B_K$ such that $\bar w(\icol n]=\bar v(\icol i_n]g_n$. As $\bar v=\bar w$ we can take $i_0=0$, $i_{\ell(w)}=\ell(v)$ and $g_0=g_{\ell(w)}=e$. We argue $w$ is recognized by the path going through all $(p_{i_n},g_n)$ in order (and only those), using correct type $3$ edges.
\end{itemize}
It follows that $\tilde\Lc$ is the language recognized by $\tilde M$: $\tilde \Lc$ is regular.
\end{proof}
\begin{prop}[{Compare with \cite[Proposition 4.2]{su2020formal}}] \label{sec2:su_lemma_42}
Let $G$ be a f.g.\ group, $\Ac$ a finite generating set, and $\Rc,\Sc$ regular languages evaluating to $R,S\subseteq G$ respectively. Suppose there exists $K\ge 0$ such that, for all $g\in R\cap S$, there exist $v\in\Rc$ and $w\in\Sc$ evaluating to $g$ such that $w$ asymmetrically $K$-fellow travel to $v$. Then $R\cap S$ is rational.
\end{prop}
\begin{proof}
$\tilde \Rc\cap \Sc$ is rational and evaluates to $R\cap S$.
\end{proof}

\begin{lemma}[Nielsen basis] %\footnote{\color{red} Does this need a proof ? Can't find a reference (maybe somewhere in Ilya Kapovich's work ?). It follows the same scheme as Nielsen basis for subgroups of a free group.}
Consider $H$ a group acting freely on a simplicial tree $\Tc$, and fix a vertex $p\in\Tc$. Then $H$ is free, and admits a basis $\Nc$ such that, for any reduced word $w$ over $\Nc^\pm$ and any $1\le n\le\ell(w)$, the geodesic path from $p$ to $\overline w\cdot p$ is not covered by the geodesic paths from $p$ to $\overline w(\icol n-1]\cdot p$ and from $\overline w(\icol n]\cdot p$ to $\overline w\cdot p$.
\end{lemma}

\begin{proof}
[Proof of Theorem \ref{sec2:nightmare_bis}.]
We will only prove part (a), as part (b) is similar.

Fix $\Ac\subseteq A\cup B$ a finite symmetric generating set. Take $\Rc\subseteq \Ac^*$ a language evaluating to $R$. Consider $\Tc$ the Bass-Serre tree\footnote{For reminder, $V\Tc=G/A\sqcup G/B$ and $E\Tc=G/C$. An edge $e$ is incident to a vertex $v$ if $e\subseteq v$.} for $G=A*_CB$ and $p=A$. Take $\Nc$ a Nielsen basis for $H$ as defined above, and consider $\Sc$ the language of reduced words over $\Nc^\pm$. 

We check that, for every pair $(v,w)\in\Rc\times\Sc$ evaluating to the same point, $w$ asymmetrically $K$-fellow travel with $v$, where
\[ K=\max_{h\in\Nc}\ell_\Ac(h)+\diam_\Ac(C).\]
The previous results says that, for each $1\le n\le \ell(w)$, there exists an edge $g_nC$ such that all $\bar w(\icol j]A$ with $j<n$ lie in one component of $\Tc\setminus\{g_nC\}$, and all $\bar w(\icol j]A$ with $j\ge n$ lie in the other. Looking at the Cayley graph of $G$ (w.r.t.\ $\Ac$), this means that the (coarse) path $w$ \say{crosses} the cutset $g_nC$ exactly once, between $w(\icol n-1]$ and $w(\icol n]$. On the other side, the (genuine) path $v$ also has to cross this cutset (as it evaluates at the same endpoint as $w$). Let $i_n$ be the smallest index such that $v(\icol i_n]\in g_nC$. It's easy to check that the sequence $(i_n)_{1\le n\le \ell(w)}$ is (strictly) increasing and
\[ \dist_\Ac\big(\bar w(\icol n],\bar v(\icol i_n]\big)\le \dist_\Ac\big(\bar w(\icol n-1],\bar w(\icol n]\big)+\diam_\Ac(C) \le K. \]
\begin{center}
    \begin{tikzpicture}[scale=1.8]
	
	\pgfmathsetmacro{\sq}{sqrt(3)/2}
	\clip (-3,-2.35) rectangle (3,2.32);
	\newcommand{\hex}[4]{
		\begin{scope}[shift={(#1,#2)}, scale=#3*.98]
			\draw[#4, thick, rounded corners={#3*6pt}] (1,0) -- (.5,\sq) -- (-.5,\sq) -- (-1,0) -- (-.5,-\sq) -- (.5,-\sq) -- cycle;
		\end{scope}}
	
	\hex00{1.01}{Turquoise}
	\hex0{5/3*\sq}{2/3}{Magenta}
	\hex{5/4}{-5/6*\sq}{2/3}{Magenta}
	\hex{-5/4}{-5/6*\sq}{2/3}{Magenta}
	
	\hex{5/4}{-35/18*\sq}{4/9}{Turquoise}
	\hex{-5/4}{-35/18*\sq}{4/9}{Turquoise}
	\hex{5/6}{20/9*\sq}{4/9}{Turquoise}
	\hex{-5/6}{20/9*\sq}{4/9}{Turquoise}
	\hex{25/12}{-5/18*\sq}{4/9}{Turquoise}
	\hex{-25/12}{-5/18*\sq}{4/9}{Turquoise}
	
	\hex{25/12}{25/54*\sq}{8/27}{Magenta}
	
	\fill[LimeGreen, opacity=.6, rounded corners=2pt] (-.3,\sq-.05) rectangle (.3,\sq+.05);
	\fill[LimeGreen, opacity=.6, rounded corners=2pt, rotate=120] (-.3,\sq-.05) rectangle (.3,\sq+.05);
	\fill[LimeGreen, opacity=.6, rounded corners=2pt, rotate=240] (-.3,\sq-.05) rectangle (.3,\sq+.05);
	
	\fill[LimeGreen, opacity=.6, rounded corners=2pt, rotate=0] (-1.03,-3/2*\sq-.05) rectangle (-1.48,-3/2*\sq+.05);
	\fill[LimeGreen, opacity=.6, rounded corners=2pt, rotate=120] (-1.03,-3/2*\sq-.05) rectangle (-1.48,-3/2*\sq+.05);
	\fill[LimeGreen, opacity=.6, rounded corners=2pt, rotate=240] (-1.03,-3/2*\sq-.05) rectangle (-1.48,-3/2*\sq+.05);
	\fill[LimeGreen, opacity=.6, rounded corners=2pt, rotate=0] (1.03,-3/2*\sq-.05) rectangle (1.48,-3/2*\sq+.05);
	\fill[LimeGreen, opacity=.6, rounded corners=2pt, rotate=120] (1.03,-3/2*\sq-.05) rectangle (1.48,-3/2*\sq+.05);
	\fill[LimeGreen, opacity=.6, rounded corners=2pt, rotate=240] (1.03,-3/2*\sq-.05) rectangle (1.48,-3/2*\sq+.05);
	
	\node[circle, fill=black, inner sep=1.5pt] (e) at (-.7, -3/5*\sq) {};
	\node at (-.8, -3/5*\sq-.08) {\small$e$};
	
	\draw[Mulberry, dashed, thick] (e)
		to [out=30, in=150, looseness=.9] (.7,-3/5*\sq)
		to [out=-30, in=100] (1.15,-3/2*\sq)
		to [out=-80, in=95] (1.25, -43/18*\sq);
	\draw[Mulberry, dashed, thick]  (1.25, -43/18*\sq)
		to [out=85, in=-105] (1.35,-3/2*\sq)
		to [out=75, in=-130] (1.817,-.633*\sq)
		to [out=50, in=-95] (2.13, 1/6*\sq)
		to [out=85, in=-90] (2.17, 41/54*\sq);
	
	\draw[RedOrange, thick] (e)
	to [out=50, in=-75, looseness=4] (-.2,\sq)
	to [out=105, in=-130, looseness=2] (.5,1.9)
	to [out=50, in=35, looseness=4] (.7,1.7)
	to [out=-145, in=95, looseness=1.5] (.2,\sq)
	to [out=-85, in=135, looseness=1.5] (.85, -3/10*\sq)
	to [out=-45, in=-145] (1.7,-.35)
	to [out=35, in=-100] (2.17, 41/54*\sq);
	
	\node[circle, fill=Mulberry, inner sep=1.5pt] at (1.25, -43/18*\sq+.01) {};
	\node[Mulberry] at (1.25, -43/18*\sq-.15) {\footnotesize$\bar w(\icol 1]$};
	\node[circle, fill=black, inner sep=1.5pt] at (2.17, 41/54*\sq-.01) {};
	\node at (2.19, 41/54*\sq+.13) {\footnotesize$\bar v=\bar w$};
	
	\node[circle, fill=RedOrange, inner sep=1.5pt] (vi1) at (.85, -3/10*\sq) {};
	\node[circle, fill=RedOrange, inner sep=1.5pt] (vi2) at (1.7,-.35) {}; 
	
	\node[RedOrange, inner sep=0] (vi1text) at (1.25,0) {\footnotesize$\bar v(\icol i_1]$};
	\draw[orange, thick] (vi1) -- (vi1text);
	
	\node[LimeGreen, inner sep=0] (g1text) at (0.5,-1) {\footnotesize$g_1C$};
	\draw[LimeGreen, thick] (.6,-.82*\sq) -- (g1text);
	
	\node[Mulberry] at (-.2, -.5) {$w$};
	\node[RedOrange] at (-.4, 0) {$v$};
\end{tikzpicture}
    \captionof{figure}{The Cayley graph of ${\color{Turquoise}A}*_{\color{LimeGreen}C}{\color{magenta}B}$, and two words ${\color{RedOrange}v},{\color{Mulberry}w}$.}
\end{center}
Everything is in place for Proposition \ref{sec2:su_lemma_42}, except that $\Sc\not\subset\Ac^*$. This can be fixed replacing each letter of $\Nc^\pm$ by a geodesic representative over $\Ac$. (We get a new language $\Sc_\Ac\subset\Ac^*$, which is regular and evaluates to $H$.) We conclude that $R\cap H$ is rational.
\end{proof}
%{\color{blue} \textbf{Remark.} Asymmetric $K$-fellow travel is weaker than the usual asynchronous $K$-fellow travel. Perhaps it would be better called "oriented language $K$-quasi-convexity".}
\section{Monsters without rational cross-sections}

Let us first recall a classical lemma for regular languages
\begin{lemma}[Pumping lemma]
Let $\Lc$ be a regular language (over $\Ac$), either
\begin{itemize}[leftmargin=7mm]
    \item $\Lc$ is finite, or
    \item $\Lc$ is infinite and there exist $v_1,v_2,w\in\Ac^*$ with $w\ne \varepsilon$ such that $v_1\{w\}^*v_2\subseteq\Lc$.
\end{itemize}
\end{lemma}
A direct corollary (proven in both \cite{Gilman1987GROUPSWA, de1990groupes}) is the following: if $\Lc$ is a rational cross-section for a group $G$, then either $G$ is finite, or $G$ contain an element $\bar w$ of infinite order. In particular, infinite torsion groups don't have rational cross-sections.

There exist other dichotomies for regular languages that can be used in a similar fashion. We first need a notion of \say{smallness} for languages:
\begin{defi}[Bounded language]
A language $\Lc \subseteq \Ac^*$ is \textbf{bounded} if there exist (not necessarily distinct) words $w_1,\ldots,w_n\in \Ac^*$ such that \vspace*{-1mm}
\[ \Lc \subseteq \{w_1\}^*\{w_2\}^*\ldots \{w_m\}^*. \]
\end{defi}
A folklore dichotomy is the following:
\begin{lemma}[See for instance \cite{Tits_for_languages}] \label{sec3:true_dich}
Let $\Lc$ be a regular language (over $\Ac$), either
\begin{itemize}[leftmargin=7mm]
    \item $\Lc$ has polynomial growth and is bounded, or
    \item $\Lc$ has exponential growth and there exist $v_1,v_2,w_1,w_2\in\Ac^*$ with $w_1w_2\ne w_2w_1$ and
    \[ v_1\{w_1,w_2\}^*v_2\subseteq \Lc.\]
\end{itemize}
\end{lemma}
The corresponding notion of \say{smallness} for groups is
\begin{defi}[Bounded generation]
A group $G$ is \textbf{boundedly generated} if there exist (not necessarily distinct) elements $g_1,g_2,\ldots,g_m\in G$ such that
\[ G = \la g_1\ra \la g_2\ra \ldots \la g_m\ra. \]
\end{defi}
A direct translation of Lemma 3.3 gives
\begin{cor}\label{sec3:dichotomy}
Let $G$ be a group having some rational cross-section $\Lc$, at least one of the following must hold true:
\begin{itemize}[leftmargin=7mm]
    \item $G$ is boundedly generated
    \item $G$ contains a free submonoid $M_2$ of rank $2$.
\end{itemize}
\end{cor}
\textbf{Remark.} Note that Corollary \ref{sec3:dichotomy} is no longer a dichotomy. For instance the Baumslag-Solitar group $BS(1,2)$ has rational cross-sections in both regime of Lemma \ref{sec3:true_dich}, and indeed it is boundedly generated \textit{and} contains free submonoids.

Our goal in the next two subsections will be to construct concrete examples of groups without either properties, hence without rational cross-sections. In both cases, one condition will be easily discarded, some work being needed to reject the other.

\subsection{Extensions of infinite torsion groups}

We first look at groups mapping onto infinite torsion groups. Note that bounded generation goes to homomorphic image so, if $G$ has some quotient $T$ which is not boundedly generated (e.g.\ infinite torsion), neither is $G$. Absence of free submonoid doesn't behave as well under extension, however things can be done under conditions

\begin{thm}
Let $G$ be a group given by a short exact sequence
\[ 1 \longto N \longinto G \overset\pi\longto T \longto 1. \]
Suppose that $N\not\geqslant M_2$, and $T$ is torsion, then $G\not\geqslant M_2$. As a corollary, if $T$ is infinite, then $G$ doesn't have any rational cross-sections.
\end{thm}
\begin{proof}
Let $g_1,g_2\in G$. Let $n_1,n_2$ be the (finite) orders of $\pi(g_1)$ and $\pi(g_2)$ respectively, so that $g_1^{n_1},g_2^{n_2}\in N$. As $N$ doesn't contain any free submonoid, there exists some non-trivial \textit{positive} relation
\[ w_1(g_1^{n_1},g_2^{n_2}) = w_2(g_1^{n_1},g_2^{n_2}) \]
between $g_1^{n_1}$ and $g_2^{n_2}$ (in $N$). Obviously this relation can also be seen as a non-trivial positive relation between $g_1$ and $g_2$ (in $G$), so that no pair of elements $g_1,g_2\in G$ generates a free submonoid in $G$.
\end{proof}
\begin{cor} The following extensions don't admit any rational cross-sections
\begin{enumerate}[label={\upshape(\alph*)}, leftmargin=8mm]
    \item $\ZZ\times T$ for any infinite torsion group $T$, for instance Burnside groups $B(p,n)$ with odd exponent $p\ge 665$ and $n\ge 2$ generators, or the first Grigorchuk group $G_{(012)^\infty}$.
    
    \item $\ZZ\wr T$ for any infinite torsion group $T$. \hfill (Compare with Theorem \ref{sec5:large_torsion}.)
    
    \item Free groups in the variety $[x^p,y^p]=e$ for odd $p\ge 665$ and $n\ge 2$ generators. % {\color{red} Those are torsionfree ??}
    
    \item More generally, relatively free groups in the varieties $w_1(x_1^p,\ldots x_r^p)=w_2(x_1^p,\ldots,x_r^p)$, for distinct \textit{positive} words $w_1,w_2\in M_r$, $p\ge 665$ odd, and $n\ge 2$ generators.
\end{enumerate}
\end{cor}
\subsection{Groups acting on regular rooted trees} \label{section_monster_on_trees}

Another class of groups providing non-examples are groups with intermediate growth. Note that the growth of any cross-section gives a lower bound on the growth of the group, hence any hypothetical rational cross-section for a group of intermediate growth should have polynomial growth hence be bounded.

Known groups of intermediate growth mainly comes from two constructions: groups acting on trees following \cite{Grigorchuk_growth}, and groups \say{of dynamical origin} following \cite{Nekra_growth}. We will focus on Grigorchuk-type groups (not all of which are torsion).

Let $\Tc_d$ be the $d$-ary rooted tree, and $G\acts\Tc_d$ by automorphisms. We denote by $\Stab(\Lc_n)$ the \textit{pointwise} stabilizer of the $n$th level $\Lc_n$ of the tree, and $G_n=G/\Stab(\Lc_n)$.

\begin{prop}\label{sec3:at_most_expo}
If $G$ is boundedly generated, then $(\abs{G_n})_n$ is at most exponential.
\end{prop}
We first need a notation:
\begin{defi}
	Let $G$ be a group. We define its \textbf{exponent} as
	\[ \exp(G) = \inf\{ n >0 \mid  \forall g\in G,\; g^n=e \}\]
	(with the convention $\exp(G)=\infty$ if no such $n$ exists).
\end{defi}
\begin{proof}[Proof of Proposition \ref{sec3:at_most_expo}]
By construction $G_n$ acts faithfully on the $n$ first levels of $\Tc_d$, in particular it can be seen as a subgroup of the automorphism group of the first levels, i.e., a subgroup of the iterated \textit{permutational} wreath product
\[ W_{d,n} = \underbrace{\big(\big(\sym(d) \wr \ldots \big)\wr \sym(d)\big)\wr \sym(d)}_{n\text{ factors}} \]
It follows that elements of $G_n$ have relatively small orders, namely bounded by
\[ \exp(W_{d,n}) = \exp(\sym(d))^n. \]
Now suppose that $G$ is boundedly generated, and fix elements $g_1,g_2,\ldots,g_m\in G$ such that $G = \la g_1\ra \la g_2\ra\ldots\la g_m\ra$. Factoring by $\Stab(n)$ we get
\[ G_n = \la \overline{g_1}\ra \la \overline{g_2}\ra\ldots\la \overline{g_m}\ra\]
(where $\bar g$ is the image of $g$ in the quotient $G_n$) whence
\[ \abs{G_n} \le \abs{\la \overline{g_1}\ra}\abs{\la \overline{g_2}\ra}\ldots\abs{\la \overline{g_m}\ra} \le \exp(\sym(d))^{mn} \]
as wanted.
\end{proof}
This might seem weak, but $\abs{G_n}$ typically grows as a double exponential. This is the case as soon as the Hausdorff dimension of the closure $\overline G$ inside $W_d=\Aut(\Tc_d)$ is (strictly) positive. (Recall that the Hausdorff dimension is given by
\[ \mathrm{hdim}(\overline G) = \liminf_{n\to\infty} \frac{\log\abs{G_n}}{\log\abs{W_{d,n}}} = \liminf_{n\to\infty} \frac{\log\abs{G_n}}{d^n}\cdot\frac{d-1}{\log(d!)}\;\text{.)} \]
This includes spinal $p$-groups with $p\ge 3$ \cite{p_Spinal_dim} and weakly regular branch groups \cite{RegBranch_dim, Weakly_Regular_Branch}. For spinal $2$-groups $G_\omega$, the orders $\abs{G_n}$ have been computed in \cite{2_Spinal_dim}, and the analogous $\limsup$ is positive as soon as $\omega$ is not eventually constant, that is, as soon as $G_\omega$ is not virtually abelian. (The Hausdorff dimension itself is positive only when $\omega$ doesn't contain arbitrarily long sequences of identical symbols.)
    
We deduce that none of these groups has bounded generation. Combining this result with known results on intermediate growth among spinal groups (see \cite{Grigorchuk_growth,3_Spinal_growth}), we get new examples of groups without rational cross-sections:
\begin{cor} The following spinal groups don't admit any rational cross-section
\begin{itemize}[leftmargin=9mm]
    \item[{\upshape (a)}] Non virtually-abelian spinal groups $G_\omega$ with $p=2$. This includes the (non-torsion) Grigorchuk-Erschler group $G_{(01)^\infty}$.
    
    \item[{\upshape (b)}] Spinal groups $G_\omega$ with $p=3$ satisfying hypothesis of \cite[Theorem 4.6]{3_Spinal_growth}. This includes the Fabrykowski-Gupta group.
\end{itemize}
\end{cor}
\textbf{Remark.} The same argument works for context-free cross-sections. Indeed, just as regular languages, context-free languages are either bounded or have exponential growth \cite{Tits_for_languages}, so that groups of sub-exponential growth with context-free cross-sections should have bounded generation. This is not the case for spinal groups covered previously.

\section{Orders and wreath products: Positive results}

Let us first recall a known result that started our investigation:

\begin{prop}[See \cite{Gilman1987GROUPSWA}]
The lamplighter group 
\[ C_2\wr \ZZ = \la a,t \mid a^2= [a, t^nat^{-n}]=e \textnormal{ for } n=1,2,\ldots\ra \]
admits a rational cross-section.
\end{prop}
\begin{proof} We construct a rational cross-section over $\Ac=\{ a,t,T=t^{-1}\}$:
\[ \Lc = \Tc \sqcup \Tc a(t^+a)^*\Tc \]
where $t^+=tt^*=\{t^n:n\ge 1\}$, similarly $T^+=TT^*$, and $\Tc=t^+\sqcup T^+\sqcup \{\varepsilon\}$.
\end{proof}
Let us be a bit informal. Except for the \say{all lamps off} elements (translations) covered by the first term $\Tc$, these normal forms consist in going through the support from left to right, changing the state of each lamp to match the element represented, and never touching those lamps ever after. Similarly, for general wreath products $L\wr Q$, we would like to have some order on $Q$. Of course, this ordering should be \say{encodable} in a finite automaton, which forces the order to be left-invariant, at least under the action of some finite-index subgroup (as the only recognizable subsets of $Q$ are finite-index subgroups by a theorem of Anissimov and Seifert). This naturally leads to the following definition:
\begin{defi} \label{sec3:MLO}
A finitely generated group $G$ has \textbf{property (R+LO)} if there exists a \textit{total} left-invariant order $\prec$ on $G$ such that the associated positive cone $G_+=\{g\in G\mid g\succ e\}$ admits a rational cross-section $\Gc_+$.
\end{defi}
\textbf{Remark.} Property (R+LO) implies that $G$ admits a rational cross-section. Indeed, if $G_+$ admits a rational cross-section $\Gc_+$, then we can define
\[ \Gc_- := \{ s_n^{-1}\ldots s_2^{-1}s_1^{-1} \mid s_1s_2\ldots s_n\in\Gc_+ \} \subset (\Ac^{-1})^* \]
which is also rational (as regularity is preserved under substitution/morphism and mirror image). It follows that $G$ admits
\[ \Gc := \Gc_+\sqcup \Gc_- \sqcup \{\varepsilon\} \subset (\Ac\cup \Ac^{-1})^* \]
as a rational cross-section.

Once this (R+LO) condition defined, Gilman's construction generalizes easily.
\begin{thm} \label{sec4:positive_result}
Let $L$ and $Q$ be two groups such that
\begin{itemize}[leftmargin=5mm]
    \item $L$ has a rational cross-section;
    \item $Q$ is virtually {\upshape (R+LO)}.
\end{itemize}
Then $L\wr Q$ has a rational cross-section.
\end{thm}

\begin{proof}
Let us first suppose that $Q$ is \RLO. Let $\Qc$, $\Qc_+$ and $\Lc$ be rational cross-sections for $Q$, $Q_+$ and $L$ respectively. We define $\Lc_0=\Lc\setminus\ev^{-1}(e_L)$. We claim that
\[ \Gc := \Qc\sqcup \Qc\Lc_0(\Qc_+\Lc_0)^*\Qc \]
is a rational cross-section for $L\wr Q$:

Note that the lamplighter only moves \say{in the positive direction} from the first to the last state switch, so no switch can be undone. It follows that every word from $\Qc\Lc_0(\Qc_+\Lc_0)^*\Qc$ contains at least one non trivial state.
\begin{itemize}[label=$\star$, leftmargin=7mm]
    \item If $g\in L\wr Q$ has empty support, i.e., is a translation, then $g$ has a unique representative in $\Qc$, and no representative in $\Qc\Lc_0(\Qc_+\Lc_0)^*\Qc$.
    \item Suppose now that $g=\big((l_s)_{s\in Q},q\big)$ has non-empty support. Let $\supp g= \{s_0\prec s_1\prec \ldots\prec s_n\}$, then
    \[ g = \underset{\in \Qc}{s_0} \cdot \underset{\in \Lc_0}{l_{s_0}} \cdot \big(\underset{\in \Qc_+}{(s_0^{-1}s_1)} \cdot \underset{\in \Lc_0}{l_{s_1}}\big) \cdots \big(\underset{\in \Qc_+}{(s_{n-1}^{-1}s_n)}\underset{\in\Lc_0}{l_{s_n}}\big)\cdot \underset{\in\Qc}{(s_n^{-1}q)} \]
    is the only word of $\Gc$ mapping to $g$.
\end{itemize}
If $Q$ is only virtually (R+LO), denote by $H\leqslant Q$ a finite index subgroup (say index $n$) with property (R+LO). Note that $L^n$ admits a rational cross-section, so the first part implies that $L^n\wr H$ has a rational cross-section. Moreover $L^n\wr H$ has index $n$ is $L\wr Q$. We conclude in turn that $L\wr Q$ admits a rational cross-section.
\end{proof}

\textbf{Examples of $Q$'s satisfying the (R+LO) condition} are
\begin{enumerate}[leftmargin=7mm, label=(\alph*)]
    \item $Q=\ZZ$ with the usual order.
    \item $Q=\ZZ^d$ with lexicographic order.
    \item More generally, consider an extension $A \overset\iota\longinto B \overset\pi\longto C$. If both $A$ and $C$ have the property (R+LO), then $B$ has (R+LO) too. Indeed, the usual construction for a positive cone $B_+:=\pi^{-1}(C_+)\cup \iota(A_+)$ has a rational cross-section $\Bc_+ := \Cc_+\Ac \cup \Ac_+$.
    
    As a corollary, poly-$\ZZ$ groups (i.e., $\ZZ$-by-$\ZZ$-by-$\ldots$ groups) have (R+LO), and all finitely generated virtually nilpotent groups are virtually (R+LO).
    
    \item It is shown in \cite[Section 3.2.1]{Antoln2021RegularLO} that, if both $L$ and $Q$ admits total left-invariant \textit{rational} order, then the same is true for $G=L\wr Q$. Their argument adapts to property \RLO: If both $L$ and $Q$ have (R+LO), then $G=L\wr Q$ has (R+LO), with positive language $\Gc_+ = \Qc_+ \cup \Qc\Lc_+(\Qc_+\Lc_0)^*\Qc$.

    The positive cone is formed of elements $g=\big((l_s)_{s\in Q},q)\in L\wr Q$ such that
    \begin{itemize}
        \item either $(l_s)_{s\in Q}\equiv e_L$ and $q\in Q_+$,
        \item or $l_m\in Q_+$ where $m=\min\{ s\in Q : l_s\ne e\}$.
    \end{itemize}
    
    \item $Q=BS(1,n)=\la a,t \mid t^{-1}at=a^n\ra$ for $n\ge 1$, with positive language    \[ \Qc_+ = t^+ \sqcup \Tc a(t^+a)^* \Tc.\]
    where $\Tc=t^+\cup T^+\cup \{\varepsilon\}$. Note that other solvable Baumslag-Solitar groups $BS(1,-n)$ have index-$2$ subgroup isomorphic to $BS(1,n^2)$, so are virtually (R+LO).
    \begin{center}
        \begin{tikzpicture}

    \newcommand{\carret}[2]{
        \begin{scope}[shift={(#1,#2)}]
            \draw[thick] (0,.1) -- (0,0.28);
            \draw[thick, -latex] (0,0.62) -- (0,.9);
            \draw[-latex] (0,-.3) -- (0,-.1);
            \node at (0,0.45) {\footnotesize$t$};
            
            \ifnum 1=#2
                \draw[thick, -latex] (0,0) -- ({pow(2,-#2)},0);
            \else
                \draw[thick] (0,0) -- ({pow(2,-#2-1)-.19},0);
                \draw[thick, -latex] ({pow(2,-#2-1)+.09},0) -- ({pow(2,-#2)},0);
                \node at ({pow(2,{-#2-1})-.05},0) {\footnotesize$a$};
            \fi
        \end{scope}}
        
    \clip (-5.9,-1.85) rectangle (6,2.15);
    
    \draw[draw=red, fill=red!15, thick, rounded corners] (-.25,3) -- (-.25,.25) -- (.25,.25) -- (.25,-3.5) -- (6.2,-3.5) -- (6.2,3) -- (-.25,3);
    
\draw[-latex, thick] (-8,-2) -- (-4,-2);
\carret{-4}{-2}
\carret0{-2}
\carret4{-2}

\draw[-latex, thick] (-6,-1) -- (-4,-1);
\carret{-4}{-1}
\carret{-2}{-1}
\carret0{-1}
\carret2{-1}
\carret4{-1}

\draw[-latex, thick] (-6,0) -- (-5,0);
\carret{-5}0
\carret{-4}0
\carret{-3}0
\carret{-2}0
\carret{-1}0
\carret00
\carret10
\carret20
\carret30
\carret40
\carret50

\draw[-latex, thick] (-6,1) -- (-5.5,1);
\carret{-5.5}1
\carret{-5}1
\carret{-4.5}1
\carret{-4}1
\carret{-3.5}1
\carret{-3}1
\carret{-2.5}1
\carret{-2}1
\carret{-1.5}1
\carret{-1}1
\carret{-.5}1
\carret01
\carret{.5}1
\carret11
\carret{1.5}1
\carret21
\carret{2.5}1
\carret31
\carret{3.5}1
\carret41
\carret{4.5}1
\carret51
\carret{5.5}1

\draw[draw=none, fill=Green] (0,0) circle (2pt);
\node[red] at (1.3,-1.55) {$Q_+$};
\end{tikzpicture}
        
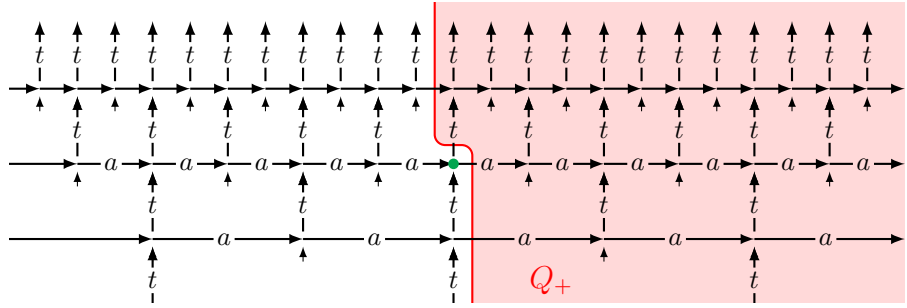
\captionof{figure}{The positive cone $Q_+=\ev(\Qc_+)$ inside $Q=BS(1,2)$}
    \end{center}
    % \item \textcolor{red}{$HNN$ extension where $\theta$ is order preserving. How should we strengthen the co-Markov condition ? $1\prec T_H\prec H_+$}
    
    \item For braid groups, the Relaxation Normal Form is regular and compatible with the Dehornoy order \cite[Theorem 5.5]{RelaxNormalForm}, hence $B_n$ is (R+LO).
    
    \newcommand{\bigast}{\scalebox{1.5}{\raisebox{-0.2ex}{$\ast$}}}
    \item A closer look at arguments of \cite{Antoln2021RegularLO} gives the following: given a finite family of groups $(G_i)$ with (R+LO), the group $\big(\bigast_i\, G_i\big)\times\ZZ$ has property (R+LO). 
\end{enumerate}

\section{Orders and wreath products: Negative results}

In this section we explore our intuition that the only way to build a rational cross-section for a wreath product $L\wr Q$ using an automaton is to switch on lamps monotonously w.r.t.\ some rational left-invariant order on $Q$. This insight translates into Lemma \ref{sec5:Positives_cones}. We deduce from this lemma a criterion (Proposition \ref{sec5:large_antichain}) ensuring some wreath products do not admit rational cross-sections, and apply it to groups similar to $C_2\wr(C_2\wr \ZZ)$. In the last subsection, we apply the criterion to wreath products $L\wr Q$ with $Q$ infinite ended.

\subsection{Main lemma}

\begin{lemma}[Positive cones on $Q$] \label{sec5:Positives_cones}
Let \vspace*{-1.5mm}
\[ 1\longto N\longinto G\overset{\pi_Q}\longto Q\longto 1 \vspace*{-1.5mm} \]
be a short exact sequence, and suppose $R\subseteq G$ has a rational cross-section $\Lc$. Let $M=(V,\Ac,\delta,*,T)$ be a \hyperlink{trimmed}{trimmed} automaton accepting $\Lc$. For each state $v\in V$, we define the language $\Lc_{v\to v}$ of words we can read along paths from $v$ to $v$. Finally, we define $P_v=\pi_Q(\ev(\Lc_{v\to v}))$. Then $P_v$ is a submonoid and
\begin{enumerate}[label={\upshape(\alph*)},leftmargin=8mm]
    \item If $N$ is torsion (i.e., $N\not\geqslant M_1$), then $P_v\cap P_v^{-1}=\{e_Q\}$ 
    
    \item If $N\not\geqslant M_2$, then $P_v\cap P_v^{-1}$ is a finite cyclic subgroup of $Q$.
\end{enumerate}
Moreover, in either case, $P_v\setminus P_v^{-1}$ is a rational sub-semigroup of $Q$ so that\vspace*{-1mm}
\[ g\prec_v h \iff g^{-1}h\in P_v\setminus P_v^{-1} \vspace*{-1mm}\]
defines a left-invariant \textbf{rational partial order} on $Q$.
\end{lemma}
\medbreak
Building a word $w\in\Lc$ correspond to following a path in the automaton $M$. Morally what our lemma says is that, as long as we stay in the strongly connected component (or \say{communication class} using Markov chain terminology) of a vertex $v\in V$, the projection of $w(\icol i]$ in $Q$ will move along a chain for $\prec_v$.

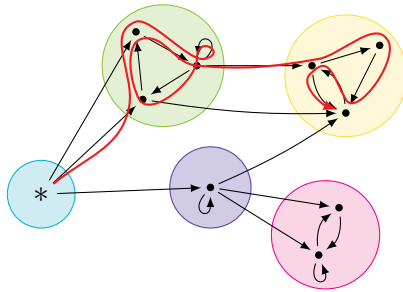
\begin{center}
    \begin{tikzpicture}[scale=.9]
    \draw[draw=Turquoise, fill=Turquoise!20] (.5,.2) circle (5mm);
    \draw[draw=Goldenrod, fill=Goldenrod!20] (5,1.95) circle (9mm);
    \draw[draw=LimeGreen, fill=LimeGreen!20] (2.3,2.1) circle (9mm);
    \draw[draw=Violet, fill=Violet!20] (3,.3) circle (6mm);
    \draw[draw=Magenta, fill=Magenta!20] (4.7,-.4) circle (8mm);
    
    \node[circle, inner sep=0pt, outer sep=2pt] (init) at (.5,.2) {$*$};
    
    \begin{scope}[every node/.style={circle, fill=black, inner sep=1pt, outer sep=2pt}]
        \node (v0) at (2,1.6) {};
        \node (v1) at (2.8,2.1) {};
        \node (v2) at (1.9,2.6) {};
     
        \node (w0) at (4.5,2.1) {};
        \node (w1) at (5.5,2.4) {};
        \node (w2) at (5,1.4) {};
    
        \node (x0) at (3,.3) {};
    
        \node (y0) at (4.9,0) {};
        \node (y1) at (4.6,-.7) {};
    \end{scope}
    
    \begin{scope}[-latex]
        \draw (init) -- (v0);
        \draw (init) -- (v2);
        \draw (v0) -- (v2);
        \draw (v2) -- (v1);
        \draw (v1) -- (v0);
        
        \draw (v1) -- (w0);
        \draw[bend right=5] (v0) to (w2);
        \draw (w0) -- (w1);
        \draw (w1) -- (w2);
        
        \draw (init) -- (x0);
        \draw (x0) -- (y0);
        \draw (x0) -- (y1);
        \draw[bend right=5] (x0) to (w2);
        
        \draw[bend right] (w0) to (w2);
        \draw[bend right] (w2) to (w0);
        
        \draw[bend left] (y0) to (y1);
        \draw[bend left] (y1) to (y0);
        
        \draw[in=85, out=40, loop, looseness=10] (v1) to (v1);
        \draw[out=-130, in=-90, loop, looseness=10] (x0) to (x0);
        \draw[out=-110, in=-70, loop, looseness=10] (y1) to (y1);
    \end{scope}
    
    \draw[Red, thick, -latex] (init)
	to [out=40, in=-100, looseness=1] (1.8,1.5)
	to [out=80, in=-160, looseness=1] (1.9,2.75)
	to [out=20, in=-170, looseness=.5] (2.9,2.2)
	to [out=10, in=90, looseness=10] (2.8,2.2)
	to [out=-90, in=-60, looseness=.8] (1.95,1.6)
	to [out=120, in=-170, looseness=.8] (2,2.5)
	to [out=10, in=175, looseness=1] (2.8,2.1)
	to [out=-5, in=-170, looseness=1] (4.6,2.2)
	to [out=10, in=125, looseness=.8] (5.6,2.5)
	to [out=-45, in=-60, looseness=.8] (5,1.6)
	to [out=120, in=45, looseness=1.5] (4.5,2)
	to [out=-145, in=160] (w2);
    \end{tikzpicture}
%\captionof{figure}{An automaton and its strongly connected components}
\end{center}
The main idea is that, as soon as $P_v\cap P_v^{-1}$ (or rather $\Tc$ in what follows) is big enough, the language will recognize a bunch of words projecting onto the same element of $q\in Q$, and embedding all this mess (aka a free submonoid) inside a single lateral class $qN$ isn't possible unless $N$ itself contains a free submonoid. (The \say{trimmed} condition is there to make sure no part of the automaton remains invisible in the recognized language.)

\textbf{Examples.}
\begin{itemize}[leftmargin=7mm]
    \item Let us first look at the rational cross-section for $F_2$ formed by all reduced words on $a,b$ and their inverses $A,B$. The group $F_2$ can be seen as a trivial extension $\{e\}\into F_2\onto F_2$, we have $N=\{e\}$ torsion, and (for the automaton in Figure \ref{auto:F2})
    \[ P_v = \{e\}\cup \{g\in F_2 \text{ whose reduced word has the form $s_1\ldots s_{\ell-1}a$ with $s_1\ne A$}\}  \]
    which does satisfy $P_v\cap P_v^{-1}=\{e\}$.
    
\begin{minipage}{.45\linewidth}
\centering
\begin{tikzpicture}[scale=1.5]
    \clip (-2,-1.5) rectangle (2,1.5);
    \node[circle, draw=Green, ultra thick, inner sep=1pt] (init) at (0,0) {$*$};
    
    \begin{scope}[every node/.style={circle, fill=Green, inner sep=2pt, outer sep= 3pt}]
        \node (a) at (1,0) {};
        \node (b) at (0,1) {};
        \node (c) at (-1,0) {};
        \node (d) at (0,-1) {};
    \end{scope}
    
    \node[circle, draw=red, thick] at (a) {};
    \node at (1.07,-.21) {\color{red}$v$};
    
    \begin{scope}[every node/.style={circle, fill=white, inner sep=1pt}]
		    \draw[-latex] (init) to node[pos=.45, inner sep=2.5pt]{$A$} (c);
		    \draw[-latex, bend right=25] (b) to node[pos=.5, inner sep=2.5pt]{$A$} (c);
		    \draw[-latex, bend left=25] (d) to node[pos=.5, inner sep=2.5pt]{$A$} (c);
		    \draw[-latex, out=210, in=150, looseness=9] (c) to node[pos=.35]{\small$A$} (c);
		    
		    \draw[-latex] (init) to node[pos=.45, inner sep=2.5pt]{$B$} (d);
		    \draw[-latex, bend right=15] (a) to node[pos=.5, inner sep=2.5pt]{$B$} (d);
		    \draw[-latex, bend left=15] (c) to node[pos=.5, inner sep=2.5pt]{$B$} (d);
		    \draw[-latex, out=300, in=240, looseness=9] (d) to node[pos=.35]{\small$B$} (d);
	
		    \draw[-latex] (init) to node[pos=.45, inner sep=2.5pt]{$a$} (a);
		    \draw[-latex, bend left=25] (b) to node[pos=.5, inner sep=2.5pt]{$a$} (a);
		    \draw[-latex, bend right=25] (d) to node[pos=.5, inner sep=2.5pt]{$a$} (a);
		    \draw[-latex, out=30, in=-30, looseness=9] (a) to node[pos=.35]{$a$} (a);
		    
		    \draw[-latex] (init) to node[pos=.45, inner sep=2.5pt]{$b$} (b);
		    \draw[-latex, bend left=15] (a) to node[pos=.5, inner sep=2.5pt]{$b$} (b);
		    \draw[-latex, bend right=15] (c) to node[pos=.5, inner sep=2.5pt]{$b$} (b);
		    \draw[-latex, out=120, in=60, looseness=9] (b) to node[pos=.35]{$b$} (b);
	\end{scope}
\end{tikzpicture}
\captionof{figure}{Automaton for $F_2=\la a,b\ra$} \label{auto:F2}
\end{minipage}
\begin{minipage}{.45\linewidth}
\centering
\begin{tikzpicture}[scale=1.5]
    \clip (-2,-.5) rectangle (2,2.5);
    \node[circle, draw=Green, ultra thick, inner sep=1pt] (init) at (0,0) {$*$};
    
    \begin{scope}[every node/.style={circle, fill=Green, inner sep=2pt, outer sep= 3pt}]
        \node (a) at (0,.9) {};
        \node (b) at (0,1.8) {};
        \node (c) at (-1,1) {};
        \node (d) at (1,1) {};
    \end{scope}
    
    \node[circle, fill=black, inner sep=1.5pt, outer sep= 3pt] (c2) at (-1.8,1.8) {};
    \node[circle, fill=black, inner sep=1.5pt, outer sep= 3pt] (d2) at (1.8,1.8) {};
    
    \node[circle, draw=red, thick] at (d) {};
    \node at (1.07,.79) {\color{red}$v$};
    
    \begin{scope}[every node/.style={circle, fill=white, inner sep=1pt}]
		    \draw[-latex] (init) to node[pos=.45, inner sep=2pt]{$a$} (a);
		    \draw[-latex] (a) to node[pos=.45, inner sep=2pt]{$a$} (b);
		    
		    \draw[-latex] (init) to node[pos=.45, inner sep=2pt]{\small$T$} (c);
		    \draw[-latex] (init) to node[pos=.45, inner sep=2pt]{$t$} (d);
		    
		    \draw[-latex] (a) to node[pos=.45, inner sep=2pt]{\small$T$} (c);
		    \draw[-latex] (a) to node[pos=.45, inner sep=2pt]{$t$} (d);
		    
		    \draw[-latex] (b) to node[pos=.5, inner sep=2pt]{\small$T$} (c);
		    \draw[-latex] (b) to node[pos=.5, inner sep=2pt]{$t$} (d);

		    \draw[-latex, bend right=20] (c) to node[pos=.5, inner sep=2pt]{$a$} (c2);
		    \draw[-latex, bend right=20] (c2) to node[pos=.5, inner sep=2pt]{\small$T$} (c);
		    
		    \draw[-latex, bend left=20] (d) to node[pos=.5, inner sep=2pt]{$a$} (d2);
		    \draw[-latex, bend left=20] (d2) to node[pos=.5, inner sep=2pt]{\small$t$} (d);
	\end{scope}
\end{tikzpicture}
\captionof{figure}{Automaton for $\la a,t\mid a^3,[a,t]\ra$}\label{auto:C3_ZZ}
\end{minipage}

	Note that the condition $N\not\ge M_2$ is indeed necessary, as $F_2$ can alternatively be seen as the extension $[F_2,F_2]\into F_2\onto \ZZ^2$ in which case $P_v=P_v^{-1}=\ZZ^2$.
	
    \item Examples with $P_v\cap P_v^{-1}$ non-trivial are quite easy to come up with. For instance, the automaton in Figure \ref{auto:C3_ZZ} recognizes a cross-section for $C_3\times \ZZ$, which can be seen as an extension $\ZZ\into C_3\times \ZZ\onto C_3$, in which case $P_v=P_v^{-1}=C_3$.
    
    \item Finally, here is a strongly connected component of the automaton constructed in  section \S4 for $C_2\wr \ZZ^2$.
    \begin{center}
\begin{tikzpicture}[scale=2]
    \begin{scope}[every node/.style={circle, fill=black, inner sep=1.5pt, outer sep= 3pt}]
        \node[fill=Green, inner sep=2pt] (init) at (0,0) {};
        \node (x) at (1,0) {};
        \node (y) at (0.5,.8) {};
        \node (Y) at (.5,-.8) {};
    \end{scope}
    
    \begin{scope}[every node/.style={circle, fill=white, inner sep=1pt}]
		    \draw[-latex, bend left=15] (init) to node[pos=.45, inner sep=2.5pt]{$s$} (x);
		    \draw[-latex, out=30, in=-30, looseness=9] (x) to node[pos=.35]{$s$} (x);
		    
		    \draw[-latex, bend right=18] (init) to node[pos=.45, inner sep=2.5pt]{$t$} (y);
		    \draw[-latex] (x) to node[pos=.5, inner sep=2.5pt]{$t$} (y);
		    \draw[-latex, out=120, in=60, looseness=9] (y) to node[pos=.35]{$t$} (y);

		    \draw[-latex] (x) to node[pos=.45, inner sep=2.5pt]{$T$} (Y);
		    \draw[-latex, out=300, in=240, looseness=9] (Y) to node[pos=.35]{$T$} (Y);
		    
		    \draw[-latex, bend left=15] (x) to node[pos=.45, inner sep=2.5pt]{$a$} (init);
		    \draw[-latex, bend right=18] (y) to node[pos=.5, inner sep=2.5pt]{$a$} (init);
		    \draw[-latex] (Y) to node[pos=.5, inner sep=2.5pt]{$a$} (init);
	\end{scope}
	
	\node[thick, circle, draw=red] at (init) {};
	\node at (-.1,-.17) {\small\color{red}$v$};
	
	\node[thick, circle, draw=orange] at (x) {};
	\node at (1.04,-.2) {\small\color{orange}$u$};
\end{tikzpicture}
\hspace{20mm}
\begin{tikzpicture}[scale=1.04]
        \clip (-2.4,-2.4) rectangle (2.3,2.27);
        
        \draw[draw=red, fill=red!15, thick, rounded corners] (-.42,3) -- (-.42,-.45) -- (.44,-.45) -- (.44,-3) -- (2.5,-3) -- (2.5,3) -- (-.45,3);
        \draw[draw=orange, pattern=north west lines, pattern color=orange, thick] (0.46, -3) rectangle (2.5,3);
        \foreach \x in {-3,...,2}{
        \foreach \y in {-3,...,2}{
            \begin{scope}[shift={(\x,\y)}]
                \node[circle, inner sep=1pt, fill=black] at (0,0) {};
                \draw[thick] (.1,0) -- (0.3,0);
                \draw[thick, -latex] (0.6,0) -- (.9,0);
                \draw[thick] (0,.1) -- (0,0.28);
                \draw[thick, -latex] (0,0.62) -- (0,.9);
                
                \node at (0.45,0) {\footnotesize$s$};
                \node at (0,0.45) {\footnotesize$t$};
            \end{scope}}}
        \draw[fill=Green, draw=Green] (0,0) circle (2.5pt);
    \end{tikzpicture}
\captionof{figure}{Component of an automaton for $C_2\wr\ZZ^2=\la a\ra\wr \la s,t\ra$ and some corresponding cones.}
\end{center}
It is notable that, even though cones obtained for different choices of $v$ are quite close in some sense (for instance, we always have
\[ w_{u\to v}\,P_v\,w_{v\to u} \subseteq P_u \subseteq w_{v\to u}^{-1}\, P_v \,w_{u\to v}^{-1} \]
for some $w_{u\to v},w_{v\to u}\in Q$), they can have drastically different behavior ($P_v$ gives rise to a total order $\prec_v$ while $\prec_u$ has \hyperlink{sec5:cd}{\textit{chain density}} $\to 0$ w.r.t.\ balls in $\ZZ^2$).
\end{itemize}
\bigbreak
\begin{proof}[Proof of Lemma \ref{sec5:Positives_cones}.] Let us define 
    \[ \Tc = \{ w \in \Lc_{v\to v} \mid \pi_Q(\bar w)\in P_v\cap P_v^{-1} \}.\]
Both cases are pretty similar
\begin{enumerate}[label=(\alph*), leftmargin=7mm]
    \item Let $w_+\in \Tc$. Denote $g=\pi_Q(\bar w_+)\in P_v\cap P_v^{-1}$, let $w_-\in\Tc$ satisfying $\pi_Q(\bar w_-)=g^{-1}$ and define $w=w_+w_-$. Note that $w\in \Lc_{v\to v}$ and $\bar w\in N$. Since the automaton is trimmed there exist paths from $*$ to $v$, and from $v$ to a terminal state. Let $v_0,v_1\in\Ac^*$ be words we can read along such paths, we get
    \[ v_0\{w\}^*v_1\subseteq \Lc.\]
    If $w\ne\varepsilon$ the language on the right hand side is infinite. As $\ev\colon\Lc\to G$ is injective, we get an infinite order element $\bar w$ in $N$, absurd! So the only possibility is $w=w_+=w_-=\varepsilon$: we get $\Tc=\{\varepsilon\}$ and $P_v\cap P_v^{-1}=\pi_Q(\ev(\Tc))=\{e_Q\}$.
    
    \newpage
    
    \item We first show that $P_v\cap P_v^{-1}$ is a torsion subgroup.
    
    Define $g$, $w_+$, $w_-$, $v_0$ and $v_1$ as in part (a), with the extra assumption that $g$ has infinite order. In particular, there does not exist another element $h\in Q$ such that both $g$ and $g^{-1}$ are positive powers of $h$, hence the same holds for $w_+$ and $w_-$. It follows that $w_+w_-\ne w_-w_+$ hence $\{w_+w_-,w_-w_+\}^*\leqslant \Lc_{v\to v}$ is a free monoid. Since
    \[ v_0\{w_+w_-,w_-w_+\}^*v_1\subseteq \Lc \]
    and $\ev\colon\Lc\to G$ is injective, we get a submonoid $\ev\{w_+w_-,w_-w_+\}^*\leqslant N$, absurd! We conclude that all elements in $P_v\cap P_v^{-1}=\pi_Q(\ev(\Tc))$ are torsion.
    
    Suppose there exist $w_1,w_2\in \Tc$ s.t.\ $w_1w_2\ne w_2w_1$. Let $n_i$ be the (finite) order of $g_i=\pi_Q(\bar w_i)$. We get a free monoid $\{w_1^{n_1},w_2^{n_2}\}^*\le \Lc_{v\to v}$ evaluating to a free submonoid $\ev\{w_1^{n_1},w_2^{n_2}\}^*\le N$, absurd! Hence $\Tc$ is a commutative submonoid of $\Ac^*$, i.e., there exists $w_0\in\Ac^*$ such that $\Tc\subseteq \{w_0\}^*$. It follows that $P_v\cap P_v^{-1}$ is a cyclic subgroup (generated by some power of $\pi_Q(\bar w_0)$).
\end{enumerate}
Note that, in both cases, $\Tc$ is a regular language (submonoids of $\{w_0\}^*\simeq\NN$ are finitely generated hence regular), so $\Lc_{v\to v}\setminus\Tc$ is regular too, and evaluates to $P_v\setminus P_v^{-1}$ in $Q$.
\end{proof}

\subsection{A criterion for wreath products}

Our interest will be directed towards wreath products $G=L\wr Q$, in which case $N=\bigoplus_Q L$. We first prove that $N$ contains free submonoid if and only if $L$ does.

\begin{prop}\label{sec5:direct}
Let $(G_i)_{i\in I}$ be monoids. Suppose none of the $G_i$ contains free (non-abelian) submonoids, then $\bigoplus_{i\in I}G_i$ doesn't contain free submonoids either.
\end{prop}
We begin by a general result on relations:
\begin{lemma}
Suppose two elements $g_1,g_2\in G$ do not generate a free monoid, then they satisfy a relation $w_1(g_1,g_2)=w_2(g_1,g_2)$ with $w_1, w_2\in M_2=\{x,y\}^*$ two distinct positive words such that $\ell(w_1)=\ell(w_2)$.
\end{lemma}
\begin{proof}[Proof of the lemma]
First of all, $g_1,g_2$ do not generate a free submonoid, so they satisfy a relation $v_1(g_1,g_2)=v_2(g_1,g_2)$ for some $v_1\ne v_2\in M_2$. If $\ell(v_1)=\ell(v_2)$ we're done. Suppose without lost of generality that $\ell(v_1)>\ell(v_2)$, and that $x$ is the $(\ell(v_2)+1)$-th letter of $v_1$. Let us take $w_1 = v_1yv_2$ and $w_2=v_2yv_1$. Obviously
\[ w_1(g_1,g_2) = v_1(g_1,g_2)\,g_2\,v_2(g_1,g_2) = v_2(g_1,g_2)\,g_2\,v_1(g_1,g_2) = w_2(g_1,g_2).\]
Moreover the $(\ell(v_2)+1)$-th letter of $w_1$ (resp. $w_2$) is $x$ (resp. $y$), whence $w_1\ne w_2$.
\end{proof}

\newpage

\begin{proof}[Proof of Proposition \ref{sec5:direct}.] We first deal with the case $\abs I=2$, so the direct sum can be written as $G\times H$. Fix $x_1,x_2\in G\times H$, say $x_i=(g_i,h_i)$. Let us construct a non-trivial positive relation between $x_1$ and $x_2$.
\begin{itemize}[leftmargin=6mm]
    \item As $G$ doesn't contain a free submonoid, there exist distinct words $u_1,u_2\in M_2$ with $\ell(u_1)=\ell(u_2)=:l$ and $u_1(g_1,g_2)=u_2(g_1,g_2)=:\tilde g$. Let $\tilde h_i:=u_i(h_1,h_2)$.

    \item As $H$ doesn't contain a free submonoid, there exist distinct words $v_1,v_2\in M_2$ with $\ell(v_1)=\ell(v_2)$ and $v_1(\tilde h_1,\tilde h_2)=v_2(\tilde h_1,\tilde h_2)$. 

    \item Consider $w_i(x,y) = v_i(u_1(x,y),u_2(x,y))$ for $i=1,2$. The announced relation is
    \[ w_1(x_1,x_2) = w_2(x_1,x_2) \]
    This equality clearly holds in the second component. In the first component it reads as $v_1(\tilde g,\tilde g) = v_2(\tilde g,\tilde g)$ which follows from $\ell(v_1)=\ell(v_2)$. Moreover this relation is non-trivial. Indeed, $v_1$ and $v_2$ differs on some letter, wlog the $j$-th letter is $x$ in $v_1$ and $y$ in $v_2$, then $w_i\big((j-1)l\col  jl\big]=u_i$ for $i=1,2$, but $u_1\ne u_2$ so that $w_1\ne w_2$.
\end{itemize}
The more general case where $I$ is finite comes by induction from the case $\abs I=2$. Finally the result extends to arbitrary sums using that \say{not containing a free submonoid} is a local property, hence goes to direct limits.
\end{proof}

We are now ready to prove our criterion
\begin{prop}\label{sec5:large_antichain}
Let $Q$ be a finitely generated group. Suppose that, for any finite sequence of left-invariant rational partial orders $\prec_1,\prec_2,\ldots,\prec_n$ on $Q$, there exists an arbitrarily big set $S\subseteq Q$ which is an antichain w.r.t.\ to all orders $\prec_i$. Then $L\wr Q$ doesn't have any rational cross-section for any non-trivial group $L\not\ge M_2$.
\end{prop}
Recall that a set $S\subset Q$ is an \textbf{antichain w.r.t.\ an order $\prec$} if and only if it doesn't contain $p,q\in S$ such that $p\prec q$ (i.e., distincts elements of $S$ are always incomparable).
\begin{proof}
For the sake of contradiction, let $\Lc$ be a rational cross-section for $L\wr Q$. WLOG assume $\Ac\subset L\cup Q$. Let $M=(V,\Ac,\delta,*,T)$ be a deterministic trimmed automaton accepting $\Lc$. Consider all orders $\prec_v$ for $v\in V$, and let $S\subseteq Q$ with 
\[ \abs S>\sum_{v\in V} \abs{P_v\cap P_v^{-1}} \vspace*{-2mm}\]
be a large common antichain. We consider a non-trivial element $h\in L$, and define $g=(h\cdot \mathds 1_S, e_Q)\in L\wr Q$. We show that no element in $\Lc$ represents $g$.

By contradiction, suppose that $g=\overline w$ for some $w\in\Lc$. For all $s\in S$, there exists $1\le i_s\le \ell(w)$ s.t.\ the state of the lamp on the site $s$ changes between $w(\icol i_s-1]$ and $w(\icol i_s]$. Recall that $\Ac\subset L\cup Q$ hence the lamplighter can only change the state of the lamp he is standing next to. In other words
\[ \pi_Q\big(\bar w(\icol i_s-1]\big) = \pi_Q\big(\bar w(\icol i_s]\big)=s. \]
By the pigeonhole principle, there exist $s,t\in S$ such that following both prefixes $w(\icol i_s]$ and $w(\icol i_t]$ in the automaton lead to the same state $v\in V$ and such that $s^{-1}t\not\in P_v\cap P_v^{-1}$. However
\[ w(i_s\col i_t] \in \Lc_{v\to v} \implies s^{-1}t = \pi_Q\big(\bar w(i_s\col i_t]\big) \in P_v,\]
so that $s^{-1}t\in P_v\setminus P_v^{-1}$, i.e., $s\prec_v t$, contradiction!
\end{proof}

\begin{thm}\label{sec5:large_torsion}
Suppose $Q$ has arbitrarily large or infinite torsion subgroups. Then $L\wr Q$ doesn't have any rational cross-section for any non-trivial group $L\not\ge M_2$.
\end{thm}
\begin{proof}
Distinct elements $s,t$ of a given torsion subgroup can never be comparable w.r.t.\ any left-invariant order, as $s^{-1}t$ has finite order. Put another way, large torsion subgroups form large common antichains.
\end{proof}

\bigbreak

\textbf{Remark.} It is natural to ask whether Proposition \ref{sec5:large_antichain} can be improved all the way to a genuine reciprocal of Theorem \ref{sec4:positive_result}, i.e., is the following statement true?

\begin{adjustwidth}{3mm}{3mm}
\textbf{Conjecture A:} The group $L\wr Q$ (with $L\not\simeq\{e\}$) has a rational cross-section if and only if $Q$ is virtually (R+LO) and $L$ has a rational cross-section.
\end{adjustwidth}

Getting back information on $L$ from \say{$L\wr Q$ has a rational cross-section} seems quite difficult. For instance, even for $Q=C_2$, it reduces (up to known results) to \say{$L\times L$ has a rational cross-section only if $L$ does} which is open. For this reason, we propose

\begin{adjustwidth}{3mm}{3mm}
\textbf{Conjecture B:} $C_2\wr Q$ has a rational cross-section iff $Q$ is virtually (R+LO).
\end{adjustwidth}

Indeed, getting back information on $Q$ seems more doable. A further argument toward the conjecture is the following strengthening of Proposition \ref{sec5:large_antichain}: let $\prec$ be a partial order on $Q$ and $S\Subset Q$ a finite subset. We define the \hypertarget{sec5:cd}{\textit{chain density}} of $\prec$ as
\[ CD(\prec,S) = \frac1{|S|}\max\big\{|C| : C\subseteq S\textit{ is a chain w.r.t.}\prec \big\} \]
\begin{prop}
Suppose $L\wr Q$ has a rational cross-section $\Lc$ with $L\not\ge M_2$ non-trivial. Let $M=(V,\Ac,\delta,*,T)$ be a finite automaton recognizing $\Lc$. There exists $\varepsilon=\varepsilon(M)>0$ such that, for all $S\Subset Q$, there exists $v\in V$ such that
\[ CD(\prec_v,S) \ge \varepsilon. \]
\end{prop}
Do these inequalities imply that one of these left-invariant orders restricts to a \textbf{total} order on a finite index subgroup of $Q$? This is true whenever $Q=\ZZ^d$ for instance. Some equivariant version of Dilworth's theorem might be useful.
% \rem With adapted generating sets, $C_n\wr(C_n\wr \ZZ)\to \ZZ\wr(\ZZ\wr \ZZ)$ (which is (M) by remark \thesection.6(f)) in the space of marked groups, hence property (M) isn't open.

\subsection{Rational cones in free and infinite-ended groups}

In this section, we prove the following strengthening of a result by Hermiller and Sunic from \cite{hermiller2016positive}. We need a slight adaptation of their argument to deal with partial orders.
\\
\begin{prop}\label{sec4:antichain}
Let $F_2=\la a,b\ra$ be a free group. For any rational order $\prec$ on $F_2$, there exists $S\leqslant F_2$ of rank $2$ such that $S$ is an antichain w.r.t.\ $\prec$.
\end{prop}
\begin{proof}
Let $P=\{g\in F_2 \mid g\succ e\}$. We suppose on the contrary that any subgroup $S\leqslant F_2$ of rank $2$ intersects $P$, aiming for a contradiction.

Let $\Ac=\{a,b,a^{-1},b^{-1}\}$, and let $\Pc\subseteq \Ac^*$ be a regular language for $P$. By Benois' lemma, we may assume $\Pc$ consists only of reduced words over $\Ac$. Let $\abs V$ be the number of states in a corresponding automaton.

$\star$ We show that, for any $g\in F_2$, there exists $g'\in P$ such that $\dist(g^{-1},g')<\abs V$.

Let $w\in\Ac^*$ be the reduced word for $g^{-1}$. Suppose wlog that $w$ ends with $b^{\pm1}$, hence 
\[ waba^{-1}w^{-1},\;wab^{-1}a^{-1}w^{-1},\;wa^2ba^{-2}w^{-1}\quad\text{and}\quad wa^2b^{-1}a^{-2}w^{-1} \]
are all reduced words, and $S=\la g^{-1}aba^{-1}g, g^{-1}a^2ba^{-2}g\ra$ is a rank $2$ subgroup of $F_2$. By our assumption there exists $k\in P\cap S$, with corresponding reduced word
\[ wa\ldots a^{-1}w^{-1} \in \Pc.\]
As a corollary, $w$ is a prefix of a word in $\Pc$, hence there exists another word $v\in\Ac^*$ of length $\ell(v)<\abs V$ such that $wv\in\Pc$. Finally $g'=\ev(wv)$ has the announced properties.

$\star$ The sequence defined by $g_0=e$ and $g_{n+1}=g_ng_n'$ for all $n\ge0$ is an infinite chain (w.r.t.\ $\prec$) supported in $B(e,\abs V)$ (which is finite), contradiction!
\end{proof}

\begin{thm} \label{sec5:infinite_ended}
Let $Q$ be a group with infinitely many ends. Then $L\wr Q$ doesn't have any rational cross-section for any non-trivial group $L\not\ge M_2$.
\end{thm}
\begin{proof}
We use once again Proposition \ref{sec5:large_antichain}, we just have to provide large antichains. We first reduce to the case $Q=F_2$ and then provide an infinite antichain in $F_2$.

$\star$ Let $\prec_1,\prec_2,\ldots,\prec_n$ be a finite sequence of rational orders on $Q$ with corresponding (rational) positive cones $P_1,P_2,\ldots,P_n$.

Using Stallings' classification of groups with infinitely many ends, we know $Q$ is either a free amalgamated products $A*_CB$ over a finite subgroup $C$, or an HNN extension $A*_C t$ over finite subgroups $\iota,\iota'\colon C\into A$. In either case standard Bass-Serre theory gives a free subgroup $F_2\le Q$ of rank $2$ acting freely on the corresponding Bass-Serre tree. Using Theorem \ref{sec2:nightmare_bis}, we get that all intersections $P_i\cap F_2$ are rational: all these orders restrict to rational orders on $F_2$.

$\star$ It remains to show that, given rational orders $\prec_1,\prec_2,\ldots,\prec_n$ on $F_2$, we can find an infinite antichain. We show by induction on $n$ that there exists a subgroup $S\le F_2$ of rank $2$ which is a common antichain for all those orders.

The case $n=1$ is covered by Proposition \ref{sec4:antichain}. Suppose our hypothesis holds for $n-1$, and let $T\le F_2$ of rank $2$ be a common antichain for $\prec_1,\ldots,\prec_{n-1}$. Let $P_n$ be the (rational) positive cone of $\prec_n$. Using Benois' lemma, we know that $P_n\cap T$ is rational: $\prec_n$ restricts to a rational order on $T$. Using Proposition \ref{sec4:antichain} once more, we get a subgroup $S\le T$ of rank $2$ which is an antichain for $\prec_n$, hence for all $\prec_1,\ldots,\prec_n$.
\end{proof}

\textbf{Remark.} Let us informally explain why Lemma \ref{sec5:Positives_cones} is so effective on wreath products. Given a word $w\in S^*$, we can see it as a path (starting at $e$) in the Cayley graph $\Cay(G,S)$. We can also consider the projected path in any quotient $Q$. For generic extensions, asking for a word $w$ to represent an element $g\in G$ only imposes the endpoint of the projected path (namely $\pi_Q(g)$). However, for $G=L\wr Q$, more is true: the projected path should also pass by all sites $s\in Q$ with \say{lamps on} (for a well-chosen $g$, this could be all sites in a large antichain), and this can be particularly hard to achieve when building $w$ following a path in an automaton.

%That being said, wreath products are not the only extensions where this observation holds. For instance, the following extension $G$ of $C_2*C_3$ is another group without rational cross-section. Elements of $G$ are equivalence classes of words over $\{a,b^\pm\}$, seen as paths in the following Cayley graph. Two words are equivalent if the path $w_1w_2^{-1}$ is closed and has winding number $0$ w.r.t.\ \emph{each} triangle in the graph. 

%\import{Source_files/}{tikz_wreath_like.tex}
	
%We can then go back and adapt the previous proof (with minor adaptations only needed in Proposition \ref{sec5:large_antichain}). $G$ is indeed an extension
%\[ 1 \longto \bigoplus_{(C_2*C_3)/C_3}\ZZ \longinto G \longto C_2*C_3 \longto 1 \vspace*{-2mm} \]
%with $\bigoplus \ZZ\not\ge M_2$, and $C_2*C_3$ infinite ended.

%%%%%%%%%%%%%%%%%%%%%%%%%%%%%%%%%%%%%%%%%%%%%%%%%%%%%%%%%%%

% \rem{Combined with the fact $\ZZ\wr B_3$ has property (M), this gives the only example (at least that I know of) of groups $\ZZ\wr F_2\leqslant \ZZ\wr B_3$ where the subgroup doesn't have property (R) while the larger one has property (R).}

% \rem{By residual finiteness of $F_2$, there exists a sequence of finite groups $(G_n)_n$ s.t. $G_n\to F_2$ in the space of marked groups (with adapted generating sets), hence $\ZZ\wr G_n\to \ZZ\wr F_2$. It follows property (M) isn't closed.}
\section{Torsion-by-$\ZZ$ groups and Houghton's $H_2$}

Particularizing Lemma \ref{sec5:Positives_cones} to torsion-by-$\ZZ$ groups, we get that $P_v\subseteq\ZZ_{\ge 0}$ or $P_v\subseteq\ZZ_{\le 0}$. Moreover the only word in $\Lc_{v\to v}$ evaluating to $0$ is the empty word $\varepsilon$. This means that, following a path in a strongly connected component of the automaton $M$, we should go up (resp. down) in the $\ZZ$-coordinate every so often, with only finitely many options in between those steps up (resp. down). More precisely, we have

\begin{prop} \label{sec6:torsion-by-Z}
Let $G$ be a group given by a short exact sequence
\[ 1\longto N\longinto G\overset{\pi}\longto \ZZ\longto 1 \]
with $N$ \textbf{torsion}. Suppose $R\subseteq G$ has a rational cross-section and fix any $t\in\pi^{-1}(1)$. Then there exists a finite subset $S\subseteq N$ and $m\in\NN$ such that
\begin{equation}
	R\subseteq \big((St)^*(St^{-1})^*\big)^m. \tag{$*$}
\end{equation}
\end{prop}
\textbf{Remark.} For $R$ a subgroup, condition $(*)$ is weaker than bounded generation. Indeed,
\begin{itemize}[leftmargin=7mm]
    \item If $g\in N$, then $g$ has  finite order $n$ so that
    \[ \la g\ra \subset \Big(\{e,g,g^2,\ldots,g^{n-1}\}\,t\Big)^*\Big(\{e,g,g^2,\ldots,g^{n-1}\}\,t^{-1}\Big)^* \vspace*{-1mm}\]
    \item If $g\notin N$, denote $h(g)=gt^{-\pi(g)}$. We have
    \[ \la g\ra \subseteq \Big(\{e,h(g),h(g^{-1})\}\,t\Big)^*\Big(\{e,h(g),h(g^{-1})\}\,t^{-1}\Big)^* \vspace*{-1mm} \]
\end{itemize}
It follows that bounded generation, that is the existence of $g_1,g_2,\ldots,g_m\in R$ such that $R=\la g_1\ra\ldots\la g_m\ra$, implies the existence of a finite $S\subseteq N$ such that $R\subseteq\big((St)^*(St^{-1})^*\big)^m$.
\\

\textbf{Remark.} As a reality check, let us consider $G=C_2\wr \ZZ$: the entire $R=G$ admits a rational cross-section, and can indeed be written as $R=\big((\{e,a\}t)^*(\{e,a\}t^{-1})^*\big)^2$.
\\

\begin{proof}[Proof of Proposition 6.1]
Let $\Lc$ be a rational cross-section for $R$ and $M=(V,\Ac,\delta,*,T)$ be a trimmed automaton accepting $\Lc$ (with $\Ac$ finite). Fix $m=2|V|-1$. We define
\[h\colon \begin{pmatrix} G & \longto & N \\ g & \longmapsto & gt^{-\pi(g)} \end{pmatrix}\]
Let $J:=\max_{s\in \Ac}|\pi(\bar s)|$ be the largest jump $\pi(\,\cdot\,)$ can do in one step in the automaton. For each $v\in V$, we denote its strongly connected component $K_v$, and $\Lc_{K_v\to K_v}$ the language of words we can read from any vertex in $K_v$ to any other vertex in $K_v$. Let
\[ S = \{e\}\cup \,h\!\left(\ev(\Ac) \cup \left\{ \,\bar w \;\,\Big|\begin{array}{c} w\in\Lc_{K_v\to K_v} \text{ for some }\\ v\in V\text{ and } \abs{\pi(\bar w)}\le J\end{array}\right\}\right) \]

The remainder of the proof goes as follows:
\begin{enumerate}[leftmargin=8mm, label=(\alph*)]
    \item We prove that $S$ is finite, through an upper bound on the length of possible $w$'s.
    \item We show that $\ev(\Lc_{v\to v})\subseteq (St)^*(St^{-1})^*$.
    \item We quickly check that $\ev(\Ac)\subseteq (St)^*(St^{-1})^*$.
    \item We decompose each word recognized by $M$ into a product of at most $m$ words of the previous two types, hence proving
    \[ R = \ev(\Lc) \subseteq \big((St)^*(St^{-1})^*\big)^m. \]
\end{enumerate}

(a) Consider a strongly connected component $K$, and $w\in\Lc_{K\to K}$ satisfying $\abs{\pi(\bar w)}\le J$. Suppose w.l.o.g.\ $P_v\subseteq\ZZ_{\ge 0}$ (for any - hence all - $v\in K$). Using some Loop-Erasure algorithm, we decompose any path recognizing $w$ as an union of a simple path labeled $w_0$, together with a bunch of (non-empty, simple) loops labeled $u_1,\ldots,u_r$.
\begin{center}
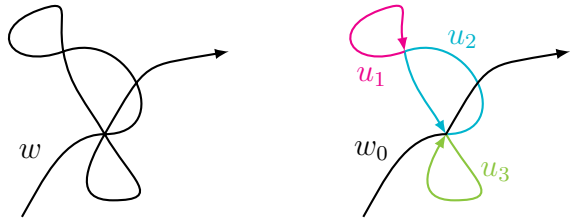

    \begin{tikzpicture}[scale=1.1]
	\draw[black, thick, -latex] (-1,-1)
		to [out=60, in=-180, looseness=.9] (0,0)
		to [out=0, in=20, looseness=2] (-.5,1)
		to [out=-160, in=-160, looseness=3] (-.8,1.5)
		to [out=20, in=100, looseness=2] (-.5,1)
		to [out=-80, in=120] (0,0)
		to [out=-60, in=0, looseness=2] (0.1,-.8)
		to [out=-180, in=-120, looseness=1.5] (0,0)
		to [out=60, in=-170, looseness=1.5] (1.5,1);
	\node[black] at (-.9,-.2) {$w$};
\end{tikzpicture}
\hspace{10mm}
\begin{tikzpicture}[scale=1.1]
	\draw[Turquoise, thick, -latex] (0,0)
		to [out=0, in=20, looseness=2] (-.5,1)
		to [out=-80, in=120] (0,0);
	\draw[magenta, thick, -latex] (-.5,1)
		to [out=-160, in=-160, looseness=3] (-.8,1.5)
		to [out=20, in=100, looseness=2] (-.5,1);
	\draw[LimeGreen, thick, -latex] (0,0)
		to [out=-60, in=0, looseness=2] (0.1,-.8)
		to [out=-180, in=-120, looseness=1.5] (0,0);
	\draw[black, thick, -latex] (-1,-1)
		to [out=60, in=-180, looseness=.9] (0,0)
		to [out=60, in=-170, looseness=1.5] (1.5,1);
		
	\node[black] at (-.9,-.2) {$w_0$};
	\node[magenta] at (-.9,.7) {$u_1$};
	\node[Turquoise] at (.2,1.15) {$u_2$};
	\node[LimeGreen] at (.6,-.45) {$u_3$};
\end{tikzpicture}
    \captionsetup{margin=6mm, font=footnotesize}
    \captionof{figure}{A path inside a component $K$ of the automaton, and its decomposition. For the decomposition, follow the path. Each time you come back at an already-visited vertex, cut the simple loop formed between the two visits, and \say{forget about the loop}, then keep going.}
\end{center}
Even though we cannot reconstruct $\overline{w}$ from $\overline{w_0}$ and the $\overline{u_i}$'s, we at least have
\[ \pi(\bar w) = \pi(\bar w_0) + \sum_{i=1}^r \pi(\bar u_i). \] % , and $\ell(w_i)\le \abs{K_v}$.
Note that $\ell(w_0)\le \abs K-1$ hence $\pi(\bar w_0)\ge J(1-\abs K)$. Recall $\pi(\bar u_i)\ge 1$ by Lemma \ref{sec5:Positives_cones}. Putting everything together $J(1-|K|)+r\le \pi(\bar w)\le J$ hence $r\le J|K|$ and finally
\[ \ell(w) = \ell(w_0)+\sum_{i=1}^r \ell(u_i) < \abs K + r \abs{K} \le J\abs K^2+\abs K.\]

(b) Let $w\in \Lc_{v\to v}$. If $\bar w=e$, this is trivial. Otherwise, Lemma \ref{sec5:Positives_cones} gives $\pi(\bar w)\ne 0$, say $\pi(\bar w)\ge 1$. Recall the notation $w(i\col j]$ for the subword of $w$ consisting of all letters from the $(i+1)$th to the $j$th, included. Let $i_0=0$. We define recursively $i_j$ as the largest integer such that $u_j=w(i_{j-1}\col i_j]$ satisfies $\pi(\bar u_j)\le J$. By definition of $J$, all $u_j$ are non-empty, so we eventually have $w=u_1u_2\ldots u_r$. By maximality of $i_j$'s (or as $\pi(\bar w)\ge 1$ whenever $r=1$), we have $1\le \pi(\bar u_j)\le J$ so that $h(\bar u_j)\in S$. Finally we get
\[\bar u_j=h(\bar u_j)t \cdot t^{\pi(\bar u_j)-1}\in (St)^*\]
for all $j$, from which $\bar w\in (St)^*$ follows. Similarly, if $\pi(\bar w)\le -1$, we have $\bar w\in (St^{-1})^*$.

(c) Let $s\in \Ac$. If $\pi(\bar s)>0$, we have $\bar s=h(\bar s)t\cdot t^{\pi(\bar s)-1}\in (St)^*$. Similarly, if $\pi(\bar s)<0$, $\bar s=h(\bar s)t^{-1}\cdot t^{\pi(\bar s)+1}\in (St^{-1})^*$. If $\pi(\bar s)=0$, we have $\bar s=h(\bar s)t\cdot t^{-1}\in (St)^*(St^{-1})^*$.

(d) Let $g\in R$. As $\Lc$ is a rational cross-section for $R$, there exists $w\in \Lc$ such that $\bar w=g$. Using another Loop-Erasure algorithm, we can rewrite
\[ w = w_1s_1w_2s_2\ldots s_{n-1}w_n \]
with $w_i\in \Lc_{v_i\to v_i}$ labeling a (possibly empty) loop at $v_i$, and $s_i$ labeling an edge from $v_i$ to $v_{i+1}$, and $n\le |V|$. We've shown $\bar w_i,\bar s_i\in (St)^*(St^{-1})^*$ which concludes the proof.
\end{proof}

\begin{center}
	\begin{tikzpicture}[scale=1.6]
	\draw[black, thick, -latex] (-.35,-.5)
		to [out=60, in=0, looseness=1.5] (-.4,-.3)
		to [out=180, in=140, looseness=1.5] (-.35,-.5)
		to [out=-40, in=-75, looseness=1.5] (0,-.1)
		to [out=105, in=-75, looseness=.9] (-.15,.75)
		to [out=105, in=160, looseness=1.5] (0,1.2)
		to [out=-20, in=30, looseness=1.5] (-.15,.75)
		to [out=-150, in=170, looseness=1.5] (0,-.1)
		to [out=-10, in=-35, looseness=2] (0.17,0.37)
		to [out=145, in=90, looseness=1.5] (-.25,0.4)
		to [out=-90, in=-120, looseness=2] (0.17,0.37)
		to [out=60, in=-120, looseness=1] (0.34,0.65);
	\node[black] at (.45,-.25) {$w$};
\end{tikzpicture}
	\hspace{10mm}
\begin{tikzpicture}[scale=1.6]
	%\draw[lightgray, step=.25] (-2,-2) grid (2,2);
	%\draw (-2,-2) grid (2,2);
	\draw[black, thick, -latex] (-.35,-.5)
		to [out=-40, in=-75, looseness=1.5] (0,-.1)
		to [out=-10, in=-35, looseness=2] (0.17,0.37)
		to [out=60, in=-120, looseness=1] (0.34,0.65);
	\draw[LimeGreen, thick, -latex]	(-.35,-.5)
		to [out=60, in=0, looseness=1.5] (-.4,-.3)
		to [out=180, in=140, looseness=1.5] (-.35,-.5);
	\draw[Turquoise, thick, -latex]	(0,-.1)
		to [out=105, in=-75, looseness=.9] (-.15,.75)
		to [out=105, in=160, looseness=1.5] (0,1.2)
		to [out=-20, in=30, looseness=1.5] (-.15,.75)
		to [out=-150, in=170, looseness=1.5] (0,-.1);
	\draw[Magenta, thick, -latex] (0.17,0.37)
		to [out=145, in=90, looseness=1.5] (-.25,0.4)
		to [out=-90, in=-120, looseness=2] (0.17,0.37);
		
	\begin{footnotesize}
		\node[black] at (.9,-.25) {$s_1s_2\ldots s_{n-1}$};
		\node[LimeGreen] at (-.55,-.2) {$w_1$};
		\node[Turquoise] at (-.5,.7) {$w_2$};
		\node[magenta] at (.075,.6) {$w_3$};
	\end{footnotesize}	
\end{tikzpicture}
	\captionsetup{margin=6mm, font=footnotesize}
	
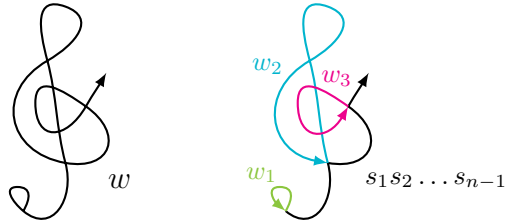
\captionof{figure}{A path in the automaton and its decomposition. For the decomposition, start at the starting vertex $*$ and skip directly to the last visit of $*$, hence bypassing a (possibly empty) loop $*\to *$, then go to the next vertex. Each time you enter a new vertex, skip directly to the last visit of said vertex (bypassing another loop), then keep going}
\end{center}

\subsection{No rational cross-section for Houghton's group $H_2$}

Houghton's groups form a family of groups with many interesting properties. The first member is the group $H_1 = \FSym(\NN)$ of finitely supported permutation of $\NN$, which is not finitely generated. The second is defined as
\[ H_2 = \left\{ \sigma \in \Sym(\ZZ) \;\Big|  \begin{array}{c} \exists \pi\in \ZZ \text{ such that }\sigma(x)=x+\pi \\ \text{ for all but finitely many }x \in \ZZ \end{array} \right\}. \]
It is finitely generated but not finitely presented. Higher groups in the family are finitely presented. Brown proved that $H_n$ has property $FP_{n-1}$, but not $FP_n$ (see \cite{BROWN1987}). In particular, $H_n$'s are examples of groups without finite complete rewriting system, and therefore good candidates not to have any rational cross-section. We show that $H_2$ doesn't have any rational cross-section. As a byproduct, it is not boundedly generated.
\medbreak
First note that $H_2$ is indeed a torsion-by-$\ZZ$ group, with the short exact sequence
\[ 1 \longto \FSym(\ZZ) \longinto H_2 \overset{\pi}\longto \ZZ\longto 1. \]
We show that subsets of $H_2$ of the form $\big((St)^*(St^{-1})^*\big)^m$
are \say{nicer} than $H_2$ as a whole, hence $H_2$ cannot be of this form. In order to formalize this idea we define a notion of \textit{complexity} for elements of $H_2$:

\begin{defi}
Let $h\in \FSym(\ZZ)$. We define its \textbf{crossing number} as
\[ c(h) = \max_{p\in\ZZ+\frac12} \# \{ x\in\ZZ : x<p<h(x) \} \vspace*{-2mm} \]
More generally, if $g\in H_2$, we define $c(g)=c(g t^{-\pi(g)})$.
\end{defi}

\begin{center}
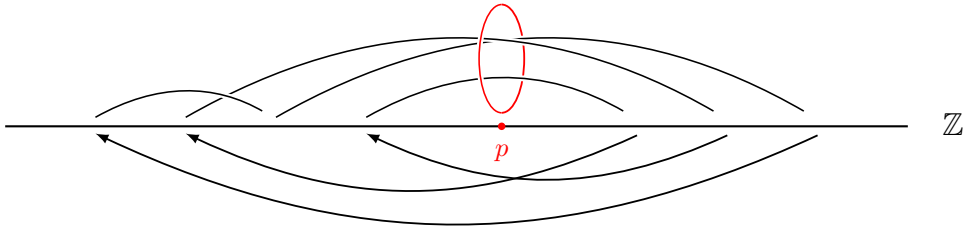

    \begin{tikzpicture}[scale=1.2, rotate=180]
        \draw[thick] (-5,0) -- (5,0);
        \node at (-5.5,-0.03) {$\ZZ$};
        
        \draw[line width=.7pt, bend left=25, -latex] (-2,0.1) to (3,0.08);
        \draw[line width=.7pt, bend left=25, -latex] (-3,0.1) to (1,0.08);
        \draw[line width=.7pt, bend left=25, -latex] (-4,0.1) to (4,0.08);
        
        \begin{scope}[thick, bend left]
        \draw[white, double=black, arrows={-latex[black]}] (1,-0.1) to (-2,-0.08);
        \draw[white, double=black, arrows={-latex[black]}] (2,-0.1) to (-4,-0.08);
        \draw[white, double=black, arrows={-latex[black]}] (3,-0.1) to (-3,-0.08);
        \draw[white, double=black, arrows={-latex[black]}] (4,-0.1) to (2,-0.08);
        \end{scope}
        
        \draw[thick, white, double=red] (-0.5,-.15) arc (90:-90:.25 and .6);
        \begin{pgfonlayer}{background}
            \draw[line width=.7pt, red] (-0.5,-.15) arc (90:270:.25 and .6);
        \end{pgfonlayer}
        
        \node[fill=red, circle, inner sep=1pt] at (-.5,0) {};
        \node[red] at (-.5,.3) {\footnotesize$p$};
    \end{tikzpicture}
    \captionof{figure}{A permutation $h\in\FSym(\ZZ)$ and a $p$ reaching the bound $c(h)=3$} \vspace*{3mm}
\end{center}
\begin{lemma}\label{sec6:complexity} Crossing numbers satisfy several properties:
\begin{enumerate}[leftmargin=8mm, label={\upshape(\alph*)}]
    \item $c(g)= c(t^mgt^n)$ for all $m,n\in\ZZ$.
    \item $c(g)=c(g^{-1})$
    \item $c(g_1g_2)\le c(g_1)+c(g_2)$
    \item $c(g)\le b-a$ for all $g\in \big(\Sym[a,b]\;t\big)^*$.
\end{enumerate}
\end{lemma}

\newpage

\begin{proof} Let us take a deep breath, and prove them one by one:
\begin{enumerate}[leftmargin=8mm, label=(\alph*)]
    \item $c(g)=c(gt^n)$ is clear (for all $g\in H_2$). Moreover, for all $h\in\FSym(\ZZ)$, we have $c(t^mht^{-m})=c(h)$, as those are \say{translated} permutations. It follows that
    \[ c(t^mgt^n) = c(t^m\,gt^{-\pi(g)}\,t^{-m}) = c(gt^{-\pi(g)})=c(g). \]
    
    \item $c(h)=c(h^{-1})$ is clear for $h\in\FSym(\ZZ)$, indeed we have
        \[ \# \{ x\in\ZZ : h(x)<p<x \} =  \# \{ x\in\ZZ : x<p<h(x) \} \]
    due to some \say{conservation of mass}. For generic $g\in H_2$, we have
    \[ c(g^{-1}) = c(g^{-1}t^{\pi(g)}) =c(t^{-\pi(g)}g) \overset{(a)}= c(gt^{-\pi(g)})=c(g) .\]
    
    \item This is clear for $h_1,h_2\in \FSym(\ZZ)$, as $x<p<h_1h_2(x)$ implies that either $h_2(x)<p<h_1h_2(x)$ or $x<p<h_2(x)$. Now for $g_1,g_2\in H_2$ we have
    \begin{align*}
        c(g_1g_2)
        & =c(g_1g_2 t^{-\pi(g_1g_2)}) = c(g_1t^{-\pi(g_1)}\; t^{\pi(g_1)}g_2 t^{-\pi(g_1g_2)}) \\
        & \le c(g_1t^{-\pi(g_1)}) + c(t^{\pi(g_1)}g_2 t^{-\pi(g_1g_2)}) \\
        & = c(g_1)+c(g_2)
    \end{align*}
    \item Let $g\in (\Sym[a,b]\;t)^n$. The associated $h=gt^{-n}\in\FSym(\ZZ)$ can be written as
    \[ h = \sigma_0 \cdot t\sigma_1t^{-1}\cdot t^2\sigma_2 t^{-2}\cdot \ldots\cdot t^{n-1}\sigma_{n-1} t^{-n+1} \]
    for some $\sigma_0,\sigma_1,\sigma_2,\ldots,\sigma_{n-1}\in\Sym[a,b]$. Note that $t^i\sigma_i t^{-i}\in \FSym(\ZZ)$ and satisfies $\supp(t^i\sigma_it^{-i})\subseteq [a+i,b+i]$. The situation can be illustrated as follows:
        \begin{center}
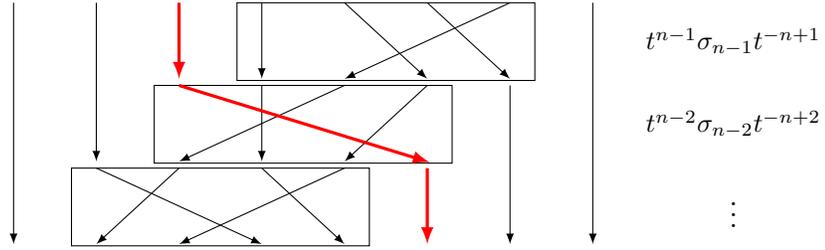

    \begin{tikzpicture}[scale=1.1]
        \newcommand{\permbox}[5]{
        \begin{scope}[shift={({-#5},{-#5})}]
            \draw (.7,0.03) rectangle (4.3,.97);
            \draw[-latex] (1,.97) -- (#1,0.05);
            \draw[-latex] (2,.97) -- (#2,0.05);
            \draw[-latex] (3,.97) -- (#3,0.05);
            \draw[-latex] (4,.97) -- (#4,0.05);
        \end{scope}}
        
        \clip (-5,-2.2) rectangle (8,1.2);
        \permbox13420
        \permbox42131
        \permbox31422
        \draw[-latex, very thick, red] (0,.97) -- (0,.05);
        \draw[-latex, very thick, red] (0,-.03) -- (3.03,-.96);
        \draw[-latex, very thick, red] (3,-1.03) -- (3,-1.95);
        
        \draw[-latex] (-2,.97) -- (-2,-1.95);
        \draw[-latex] (-1,.97) -- (-1,-.95);
        \draw[-latex] (4,-.03) -- (4,-1.95);
        \draw[-latex] (5,.97) -- (5,-1.95);
        
        \node at (6.7,0.5) {\footnotesize$t^{n-1}\sigma_{n-1}t^{-n+1}$};
        \node at (6.7,-0.5) {\footnotesize$t^{n-2}\sigma_{n-2}t^{-n+2}$};
        \node at (6.7,-1.5) {$\vdots$};
    \end{tikzpicture}
    \captionof{figure}{A \say{braid} diagram for (d)}
\end{center}
    Now observe that $h(x)\le x+(b-a)$, for all $x\in\ZZ$, which concludes. \hfill\qedhere
\end{enumerate}
\end{proof}

\newpage

Now that everything is in place, we can proceed and prove this section's main result
\begin{thm} \label{sec6:Houghton_is_down}
$H_2$ does not admit any rational cross-section.
\end{thm}
\begin{proof}
Consider a subset $R\subseteq H_2$ admitting a rational cross-section, so that $R\subseteq\big((St)^*(St^{-1})^*\big)^m$ for some finite $S\subseteq \FSym(\ZZ)$ by the previous theorem. Fix $[a,b]$ a finite interval containing the support of each $s\in S$. It follows that
\[ R \subseteq \big((\Sym[a,b]t)^*(\Sym[a,b]t^{-1})^*\big)^m \]
hence $c(g)$ is uniformly bounded by $2m(b-a)$ on $R$ (using Lemma \ref{sec6:complexity} repetitively).

On the other side, crossing numbers are not uniformly bounded on $H_2$. For instance
\[ h_K = (1,-1)(2,-2)\ldots (K,-K) \in \FSym(\ZZ)\le H_2 \]
satisfies $c(h_K)=K$ for any $K\in\NN$. The conclusion follows.
\end{proof}

\textbf{Remark.} In the previous section we used that, in order for a word $w\in \Ac^*$ to evaluate to a well-chosen $g\in L\wr Q$, the corresponding path in $Q$ should pass by a large set of elements in $Q$, which can be complicated. In this section, we show the path corresponding to any word $w$ evaluating to well-chosen $g\in H_2$ should do many back-and-forths in $Q=\ZZ$, which turns out to be just as complicated.

% {\color{blue}\textbf{Remark.} The weaker conclusion $H_2$ cannot be written as $(St)^*(St^{-1})^*$ can also be seen using the emptiness of the first BNS-invariant $\Sigma^1(H_2)$ (kinda trivial, but see Zaremsky), so that each half the Cayley graph $\Gc_{\pm\chi}=\{g\in G\mid \pm\chi(g)\ge 0\}$ with $\chi:G\to \RR:t\mapsto 1$ contains infinitely many connected components. This doesn't conclude directly for higher $m$ (and shouldn't, see examples like $G=C_2\wr \ZZ$), as we might be able to pick a connected component (among infinitely many) when diving back into $G_{-\chi}$. That being said, a finer study of usual results around $\Sigma^1$ (Morse functions ?) might be interesting, bringing a new point of view. For instance, is anything known on how those left and right connected components meet at $0$ ?}

%\textbf{Remark.} The following scheme is quite common : find some notion of \textit{complexity} of elements that's really well-behaved under exponentiation (d), and well-behaved under products (c), and then conclude it should be well-behaved over boundedly generated subsets. For instance, this is used by \cite[Corollary 3.10]{grigorchuk_1994} to prove the second bounded cohomology of boundedly generated groups is finite dimensional.

\section{A finitely presented extension of Grigorchuk's group}

The goal of this section is to prove that the finitely presented HNN-extension $\Lk$ of the first Grigorchuk group defined in \cite{Grigorchuk_1998} doesn't have a rational cross-section. This group was introduced as the first example of finitely presented group which is amenable but not elementary-amenable. In the first subsection, we recall the construction of $\Lk$, and exhibit an action of $\Lk$ on the (unrooted) $3$-regular tree. The actual proof that $\Lk$ does not admit rational cross-section, using Proposition \ref{sec6:torsion-by-Z}, comes in \S\thesection.2.

\subsection{An action of $\Lk$ on the $3$-regular tree}

Let us first recall the definition of the first Grigorchuk's group.
\begin{defi}[Grigorchuk's group, via its action by automorphisms on the infinite rooted binary tree $\{0,1\}^*$]
Grigorchuk's group is $\mathfrak G=G_{(012)^\infty}=\la a,b,c,d\ra$, where $a,b,c,d\in\Aut(\{0,1\}^*)$ are defined recursively by
\[ \begin{matrix}
\begin{cases} a(0w)=1w \\ a(1w)=0w \end{cases} \hspace*{5mm}\vspace*{2mm} &
\begin{cases} b(0w) = 0a(w) \\ b(1w)=1c(w) \end{cases} \\
\begin{cases} c(0w) = 0a(w) \\ c(1w)=1d(w) \end{cases} &
\begin{cases} d(0w) = 0w \\ d(1w)=1b(w) \end{cases}
\end{matrix} \]
\end{defi}

\newpage

We will also need the notion of section of an automorphism at a vertex:
\begin{defi}
Let $g\in\Aut(\{0,1\}^*)$ and $v\in\{0,1\}^*$. The \textbf{section of $g$ at $v$} is the unique element $g_v\in \Aut(\{0,1\}^*)$ such that
\[ \forall w\in\{0,1\}^*,\; g(vw)=g(v)g_v(w). \]
\end{defi}
For instance, the previous definition gives $a_0=\id$, $b_0=a$, $b_{10}=a$, $b_{110}=\id$, $c_1=d$. This will usually be depicted as follows:
\begin{center}
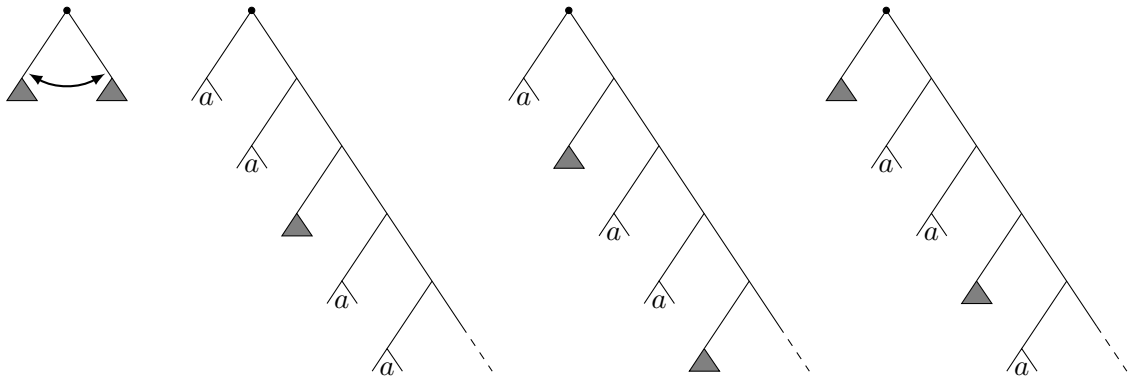

    \begin{tikzpicture}[scale=.6]
        \clip (-1.4,-8.25) rectangle (1.4,.2);
        \draw[fill=black] (0,0) circle (2pt);
        \aswitch00
    \end{tikzpicture}
    \hspace*{5mm}
    \begin{tikzpicture}[scale=.6, xscale=-1]
        \clip (-5.4,-8.25) rectangle (1.4,.2);
        \aright00
        \aright{-1}{-1.5}
        \idright{-2}{-3}
        \aright{-3}{-4.5}
        \aright{-4}{-6}
        
        \draw (0,0) -- (-4.66,-7);
        \draw[dashed] (-4.66,-7) -- (-5.33,-8);
        \draw[fill=black] (0,0) circle (2pt);
    \end{tikzpicture}
    \begin{tikzpicture}[scale=.6, xscale=-1]
        \clip (-5.4,-8.25) rectangle (1.4,.2);
        \aright00
        \idright{-1}{-1.5}
        \aright{-2}{-3}
        \aright{-3}{-4.5}
        \idright{-4}{-6}
        
        \draw (0,0) -- (-4.66,-7);
        \draw[dashed] (-4.66,-7) -- (-5.33,-8);
        \draw[fill=black] (0,0) circle (2pt);
    \end{tikzpicture}
    \begin{tikzpicture}[scale=.6, xscale=-1]
        \clip (-5.4,-8.25) rectangle (1.4,.2);
        \idright00
        \aright{-1}{-1.5}
        \aright{-2}{-3}
        \idright{-3}{-4.5}
        \aright{-4}{-6}
        
        \draw (0,0) -- (-4.66,-7);
        \draw[dashed] (-4.66,-7) -- (-5.33,-8);
        \draw[fill=black] (0,0) circle (2pt);
    \end{tikzpicture}
    \captionof{figure}{From left to right $a$, $b$, $c$ and $d$. Black triangles denotes $\id$ sections.}
\end{center}
\medbreak
Grigorchuk's group admits an $L$-presentation, exhibited by Lysenok in \cite{Lysenok1985}
\[ \Gk = \la a,c,d \mid \phi^i(a^2)=\phi^i(ad)^4 = \phi^i(adacac)^4=e \;\forall i\ge 0 \ra \]
where $\phi\colon\{a,c,d\}^*\to\{a,c,d\}^*$ is defined by $\phi(a)=aca$, $\phi(c)=cd$ and $\phi(d)=c$.

The previous $\phi$ defines an injective endomorphism $\phi\colon\Gk\to \Gk$. Pictorially we get\vspace{1mm}
\begin{center}
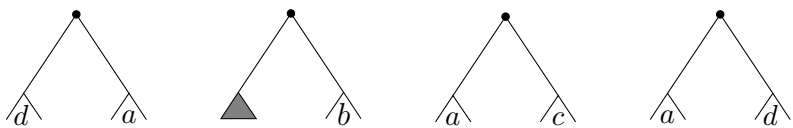

    \begin{tikzpicture}[scale=.7, xscale=-1]
        \draw[fill=black] (0,0) circle (2pt);
        \aleft00
        \dright00
    \end{tikzpicture}
    \hspace{7mm}
    \begin{tikzpicture}[scale=.7, xscale=-1]
        \draw[fill=black] (0,0) circle (2pt);
        \bleft00
        \idright00
    \end{tikzpicture}
    \hspace{7mm}
    \begin{tikzpicture}[scale=.7, xscale=-1]
        \draw[fill=black] (0,0) circle (2pt);
        \cleft00
        \aright00
    \end{tikzpicture}
    \hspace{7mm}
    \begin{tikzpicture}[scale=.7, xscale=-1]
        \draw[fill=black] (0,0) circle (2pt);
        \dleft00
        \aright00
    \end{tikzpicture}
\captionof{figure}{From left to right $\phi(a)$, $\phi(b)$, $\phi(c)$ and $\phi(d)$.} \label{phi_draw}
\end{center}
so that $\phi(g)\in\Stab(0)$ for all $g\in\Gk$, and $\phi(g)_1=g$. (This last equality proves injectivity, and will be central in order to define an action of $\Lk$.) Using this presentation, Grigorchuk \cite{Grigorchuk_1998} constructed the following finitely presented group:
\begin{defi}
The group $\Lk$ is defined as the following (ascending) HNN-extension
\[ \Lk = \Gk*_{\phi}t = \la a,b,c,d,t \;\;\Big| \begin{array}{c} a^2=b^2=c^2=d^2=bcd=(ad)^4=(adacac)^4 = e \\ t^{-1}at = aca,\, t^{-1}bt=d,\,t^{-1}ct=b,\, t^{-1}dt=c \end{array}\ra \]
\end{defi}
%\textbf{Question.} Is this group $FP_\infty$ ? Was it already known it did not admit \textbf{finite} complete rewriting system ? We show no \textbf{left-regular} complete rewriting system. \\

\newpage

\begin{prop} \label{sec7:action_on_ternary}
The group $\Lk$ acts faithfully, by automorphisms, on the $3$-regular tree.
\end{prop}
\begin{proof}[Construction] We will parametrize vertices of the $3$-regular tree by pairs $(n,w)\in \ZZ\times\{0,1\}^*$, under the extra identification $(m,1^nw)\sim (m+n,w)$. We have an edge between vertices $(n,w)$ and $(n,ws)$ for all $n\in \ZZ$, $w\in\{0,1\}^*$ and $s\in\{0,1\}$.

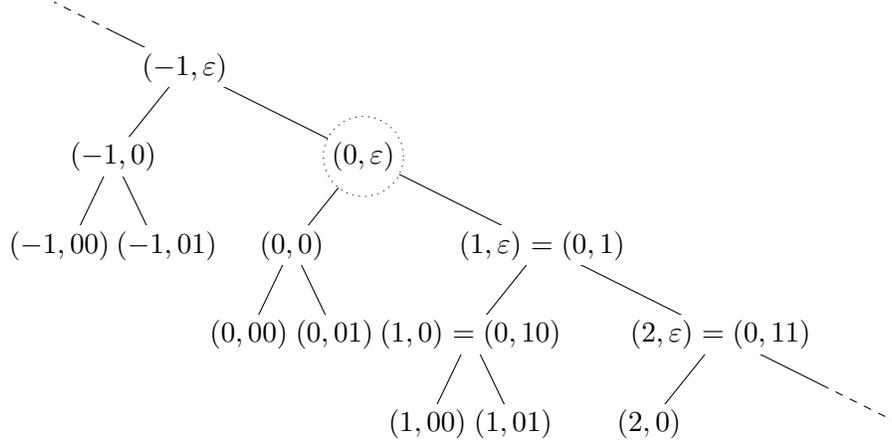
\begin{center}
    \begin{tikzpicture}[yscale=.79, xscale=-.95]
        
        \draw[dashed] (-7.5,-4.5) -- (-6.5,-3.9);
        \draw (-6.5,-3.9) -- (3.5,2.1);
        \draw[dashed] (3.5,2.1) -- (4.33,2.6);
        
        \draw (2.5,1.5) -- +(1,-1.5);
        \draw (0,0) -- +(1,-1.5);
        \draw (-2.5,-1.5) -- +(1,-1.5);
        \draw (-5,-3) -- +(1,-1.5);
        
        \draw (2.9,-1.5) -- (3.5,0) -- (4.1,-1.5);
        \draw (.4,-3) -- (1,-1.5) -- (1.6,-3);
        \draw (-2.1,-4.5) -- (-1.5,-3) -- (-.9,-4.5);
        
    \begin{scope}[every node/.style={rectangle, fill=white, inner sep=2pt}]
    \begin{small}
        \node at (-5,-3) {$(2,\varepsilon)=(0,11)$};
        \node at (-2.5,-1.5) {$(1,\varepsilon)=(0,1)$};
        \node[circle, draw, dotted, inner sep=1.5pt] at (0,0) {$(0,\varepsilon)$};
        \node at (2.5,1.5) {$(-1,\varepsilon)$};
        
        \node at (-4,-4.5) {$(2,0)$};
        \node at (-1.5,-3) {$(1,0)=(0,10)$};
        \node at (1,-1.5) {$(0,0)$};
        \node at (3.5,0) {$(-1,0)$};
        
        \node at (-2.1,-4.5) {$(1,01)$};
        \node at (.4,-3) {$(0,01)$};
        \node at (2.75,-1.5) {$(-1,01)$};
        
        \node at (-.9,-4.5) {$(1,00)$};
        \node at (1.6,-3) {$(0,00)$};
        \node at (4.25,-1.5) {$(-1,00)$};
    \end{small}
    \end{scope}
    \end{tikzpicture}
\end{center}

For $g\in\la a,b,c,d\ra$ we define $g\cdot (0,w)=(0,gw)$ through the usual action $\Gk\acts \{0,1\}^*$. Outside of this subtree, we define the action branch by branch:
\begin{itemize}[leftmargin=7mm]
    \item We start by $a\cdot (-1,0w)=(-1,0(d\cdot w))$ and $a\cdot (-2,0w)=(-2,0(dad\cdot w))$, and \\
   	$a\cdot (n,0w)=(n,0(x\cdot w))$ with $x=e,a,a$ for $n\equiv 0,1,2\pmod 3$ resp. (for $n<-2$)
    \item $b\cdot(n,0w)= (n,0(x\cdot w))$ with $x=a,a,e$ for $n\equiv 0,1,2\pmod 3$ respectively,
    \item $c\cdot(n,0w)= (n,0(x\cdot w))$ with $x=a,e,a$ for $n\equiv 0,1,2\pmod 3$ respectively,
    \item $d\cdot(n,0w)= (n,0(x\cdot w))$ with $x=e,a,a$ for $n\equiv 0,1,2\pmod 3$ respectively,
\end{itemize}
Finally $t\cdot (n,w)=(n-1,w)$.

In particular the $\Gk$-action on each branch $\{(n,0w)\mid w\in\{0,1\}^*\}$ factors through $D_8=\la a,d\ra$ for $n=-1$, through $V_4=\la a,dad\ra$ for $n=-2$, and through $C_2=\la a\ra$ for $n<-2$ (the morphism $\psi_n$ onto $C_2$ depending only on $n\bmod 3$).
\[ \begin{matrix}
    \psi_{-1}(a) = d \hspace*{4.2mm}\\
    \psi_{-2}(a) = dad  \\
    \psi_0(a)=\psi_2(b)=\psi_1(c)=\psi_0(d)=e
\end{matrix} \quad\text{and}\hspace{8mm}  {\setlength\arraycolsep{2pt} \begin{matrix}
& & \psi_1(a) & = & \psi_2(a) & = a \\
\psi_0(b) & = & \psi_1(b) & & & = a \\
\psi_0(c) & &  & = & \psi_2(c) & = a \\
& & \psi_1(d) & = & \psi_2(d) & = a
\end{matrix}} \]
This defines an action of $F(a,b,c,d,t)$, remains to check that each relation is satisfied
\begin{itemize}[leftmargin=7mm]
    \item The first line of relations holds in $\Gk$ (so holds in the subtree below $(0,\varepsilon)$), obviously $a^2, b^2, c^2, d^2$ acts as the identity on each branch, and $x^4=e$ is a law in $D_8$
    \item Relations $t^{-1}bt=d$, $t^{-1}ct=b$ and $t^{-1}dt=c$ are trivial.
    \item $t^{-1}at=aca$ is equivalent to a bunch of conditions on the actions $\psi_{-i}(a)$. First, we should have $\psi_{-1}(a)=d$ (compare with the left diagram of Figure \ref{phi_draw}). Then
    \[ \forall i\ge 1,\; \psi_{-(i+1)}(a) = \psi_{-i}(a)\psi_{-i}(c)\psi_{-i}(a) \]
    so $\psi_{-2}(a)=dad$, $\psi_{-3}(a)=dad\, e\, dad=e$, $\psi_{-4}(a)=eae=a$, and so on.
\end{itemize}
\begin{center}
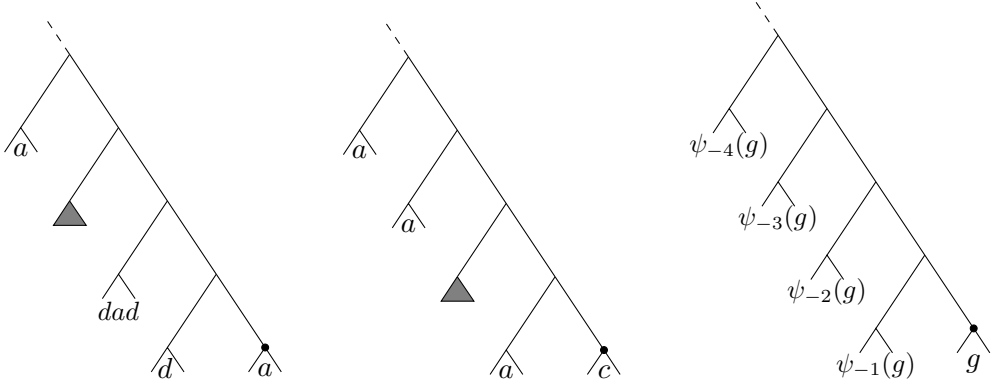

    \begin{tikzpicture}[scale=.65, xscale=-1]
        \draw[fill=black] (0,0) circle (2pt);
        \aleft1{1.5}
        \dright1{1.5}
        \emptright23
        \node at (3,.75) {\footnotesize$dad$};
        \idright3{4.5}
        \aright46
        \draw (1,1.5) -- (4,6);
        \draw[dashed] (4,6) -- (4.5,6.75);
    \end{tikzpicture}
    \hspace{5mm}
    \begin{tikzpicture}[scale=.65, xscale=-1]
        \draw[fill=black] (0,0) circle (2pt);
        \cleft1{1.5}
        \aright1{1.5}
        \idright23
        \aright3{4.5}
        \aright46
        \draw (1,1.5) -- (4,6);
        \draw[dashed] (4,6) -- (4.5,6.75);
    \end{tikzpicture}
    \hspace{5mm}
    \begin{tikzpicture}[scale=.65, xscale=-1]
        \draw[fill=black] (0,0) circle (2pt);
        \emptleft1{1.5}
        \node at (0,-.73) {\small$g$};
        \emptright1{1.5}
        \node at (2,-.75) {\footnotesize$\psi_{-1}(g)$};
        \emptright23
        \node at (3,.75) {\footnotesize$\psi_{-2}(g)$};
        \emptright3{4.5}
        \node at (4,2.25) {\footnotesize$\psi_{-3}(g)$};
        \emptright46
        \node at (5,3.75) {\footnotesize$\psi_{-4}(g)$};
        \draw (1,1.5) -- (4,6);
        \draw[dashed] (4,6) -- (4.5,6.75);
    \end{tikzpicture}
    \captionof{figure}{Action of $a$, $c$ and a random $g\in\Gk$ on the tree. The highlighted vertex is $(0,\varepsilon)$}
\end{center}
So we get a genuine $\Lk$-action. Moreover, the action restricted to $t^{n}\Gk t^{-n}$ is faithful (look under $(-n,\varepsilon)$), and any $g\notin N=\bigcup_n t^n\Gk t^{-n}$ shifts levels so acts non-trivially.
\end{proof}

\begin{rmdef}
Our tree has slightly more structure as it is graduated. Indeed we can define the \textbf{level of vertex} $(n,v)$ as $n+\ell(v)$, where $\ell(v)$ is the length of $v$. In particular we can define a relation \say{is a descendant of}, which is preserved by the action. Therefore we can still define the \textbf{section of an automorphism} $g$ at a vertex $(n,v)$ as the unique element $g_{(n,v)}\in \Aut(\{0,1\}^*)$ satisfying
\[ g\cdot (n,vw)= (\tilde n,\tilde v \,g_{(n,v)}w), \text{ where }(\tilde n,\tilde v)=g\cdot (n,v). \]
\end{rmdef}

\textbf{Remark.} The boundary of the tree can be identified with the set of doubly infinite sequences of $0$ and $1$, starting with infinitely many $1$'s (together with a globally fixed end $-\infty$, which can be considered as the root at infinity). The induced action can easily be described. For $v\in \{0,1\}^*$, $s\in\{0,1\}$ and $w\in\{0,1\}^\infty$, we have
\begin{itemize}[leftmargin=7mm]
    \item $t$ shifts to the left: $t\cdot \ldots 11v|sw=\ldots 11vs|w$ 
    \item Elements $g\in \Gk$ act on the main subtree as $g\cdot \ldots11|w= \ldots11|(gw)$ where $gw$ is defined by the usual action $\Gk\acts \{0,1\}^\infty$.
    \item Finally, $g\in\Gk$ acts on other branches as $g\cdot \ldots110v|w=\ldots110(\psi_i(g)v)|(\psi_i(g)_vw)$ where $i<0$ is the position of the first $0$.
\end{itemize}

\newpage

\textbf{Remark.} The existence of an action can also be seen abstractly. Consider a sequence of nested groups $G_n$ acting on nested sets/graphs $\Omega_n$ 
\begin{center}
    \begin{tikzcd}
        G_0 \arrow[r, hook, "\phi_0"] \arrow[d, bend left=45] & G_1 \arrow[r, hook, "\phi_1"]\arrow[d, bend left=45] & G_2 \arrow[r, hook, "\phi_2"]\arrow[d, bend left=45] & G_3 \arrow[r, hook, "\phi_3"]\arrow[d, bend left=45] & \cdots \\
        \Omega_0 \arrow[r, hook, "\iota_0"] & \Omega_1 \arrow[r, hook, "\iota_1"] & \Omega_2 \arrow[r, hook, "\iota_2"] & \Omega_3 \arrow[r, hook, "\iota_3"] & \dots
    \end{tikzcd}
\end{center}
If all of this is equivariant, in the sense $\iota_n(g\cdot \omega)=\phi_n(g)\cdot \iota_n(\omega)$ for all $g\in G_n$, $\omega\in\Omega_n$, then we get an action $\varinjlim G_n \;\acts\;  \varinjlim \Omega_n$. If moreover all $G_n$, $\Omega_n$, $\phi_n$ and $\iota_n$ are equal to fixed $G$, $\Omega$, $\phi$ and $\iota$ respectively, then we get an action
\[ G*_\phi t = \varinjlim G_n \rtimes \la t\ra \;\acts\; \varinjlim \Omega_n.\]
Here, we just take $G=\Gk$ and $\Omega=\{0,1\}^*$, $\phi=\phi$ and $\iota\colon w\mapsto 1w$. The direct limit of all $\Gk$'s is $\bigcup_n t^n\Gk t^{-n}$, while the direct limit of $\{0,1\}^*$ is our tree. (If it wasn't clear, this shows that the graph acted on is indeed a tree since direct limits of trees are trees.)

\subsection{No rational cross-section for $\Lk$}

We apply Proposition \ref{sec6:torsion-by-Z} (rather its contrapositive), hence proving the following result.

\begin{thm} \label{sec7:Lysenok_is_down}
The extension $\Lk$ doesn't have any rational cross-section.
\end{thm}
\begin{proof} Let $N=\bigcup_n t^n\Gk t^{-n}$. Fix $S\subset N$ a finite set and $m\in \NN$. We're aiming to show
\[ \Lk \ne \big((St)^*(St^{-1})^*\big)^m.\]
Up to conjugation by some power of $t$, we may suppose $S\subset\Gk$. We also add $e\in S$. The strategy is the following: for $X\subseteq N$, we consider
\[ \sec_n(X) = \{ x_v\St(\Lc_n) \mid x\in X, v\in\Lc_0 \} \subseteq \Aut(\{0,1\}^*) / \St(\Lc_n) \]
the set of sections below vertices at level $0$, but we only care about the action down until level $n$. (Recall that $\St(\Lc_n)$ is the \textit{pointwise} stabilizer of $\Lc_n$.)

In the one hand, $\sec_n(N)$ contains the usual congruence quotient $\sec_n(\Gk)=G_n$, which has size $2^{5\cdot 2^{n-3}+2}$ (as soon as $n\ge 3$, see \cite{Grigorchuk2000}). It follows $\abs{\sec_n(N)}$ grows  as a double exponential. In the other hand, we show the cardinal of
\[ \sec_n\Big( N \cap\big((St)^*(St^{-1})^*\big)^m\Big)\]
grows exponentially in $n$, so that $\big((St)^*(St^{-1})^*\big)^m$ cannot fully cover $N$, nor $\Lk$.

\textbf{Observation.} The chain rule $(xy)_v=x_{y(v)}y_v$ gives
\[ \sec_n(XY)\subseteq \sec_n(X)\sec_n(Y)\quad\text{for }X,Y\subseteq N.\]
This allows to decompose the problem into simpler parts. Indeed, each element $h\in N\cap\big((St)^*(St^{-1})^*\big)^m $ can be written as a product $h = a_0a_1\ldots a_\ell$
with 
\[ a_i\in t^{j(i)} S t^{-j(i)} \qquad\text{for some height function } j\colon [\![1,\ell]\!]\to\ZZ \]
satisfying $j(i+1)-j(i)=\pm 1$, with the sign of this difference changing at most $2m-1$ times (and some additional boundary conditions). Pictorially,
\begin{center}
    \begin{tikzpicture}[yscale=.78, xscale=.9]
    \draw[dashed] (-4.2,0) -- (5,0);
    \draw[dashed] (-4.2,-1.5) -- (5,-1.5);
    \node (L0) at (5.4,0) {$\Gk$};
    \node (Ln) at (5.9,-1.5) {$t^{-n}\Gk t^n$};
    
    \newcommand{\drop}[4]{
        \pgfmathsetmacro{\run}{2*(#2-#4)}
        \foreach \i in {0,...,\run}{
            \pgfmathsetmacro{\x}{((\run-\i)*#1 + \i*#3)/\run}
            \pgfmathsetmacro{\y}{((\run-\i)*#2 + \i*#4)/\run}
            \pgfmathsetmacro{\ytest}{2*\y}
            \ifnum 0<\ytest
                \node[circle, red, fill=red, inner sep=1.8pt] at (\x,\y) {};
                \draw[red, thick] (\x,\y) -- (#1,#2);
            \else
                \ifnum -3<\ytest
                    \node[circle, orange, fill=orange, inner sep=1.8pt] at (\x,\y) {};
                \else
                    \node[circle, green, fill=green, inner sep=1.8pt] at (\x,\y) {};
                    \draw[green, thick] (\x,\y) -- (#3,#4);
                \fi
            \fi}}
        
    \drop{-3}{1}{-3.7}0
    \drop{-3}{1}{-1.5}{-2}
    \drop0{2}{-1.5}{-2}
    \drop0{2}2{-2.5}
    \drop{3.5}{1.5}2{-2.5}
    \drop{3.5}{1.5}{3.9}{.5}
    \node[circle, draw=black, fill=white, inner sep=1pt] at (4.1,0) {};
    \end{tikzpicture}
    \captionof{figure}{$j(i)$ versus $i$, and the bunching. Here $m=3$.}
\end{center}

If we bunch together all contiguous $a_i$ with $j(i)>0$ (red pikes), and all contiguous $a_i$ with $j(i)\le -n$ (green valleys), we get elements of $U:=\bigcup_{J\ge 1}t(St)^{J-1}(St^{-1})^J$ and $t^{-n}\Gk t^n$ respectively, so that
\[ N\cap\big((St)^*(St^{-1})^*\big)^m \subseteq S \cdot U \cdot \left( S \cdot t^{-1}St \cdot \ldots \cdot t^{n-1}St^{1-n} \cdot t^{-n}\Gk t^n \cdot t^{n-1}St^{1-n} \cdot \ldots \cdot S \cdot U\right)^{m-1}. \]
We then find estimates on the number of possibles sections for each subset:
\begin{enumerate}[leftmargin=8mm, label=(\alph*)]
    \item For $n\ge 0$, we have $\abs{\sec_n(t^{-n}\Gk t^n)}\le 266$. % (We could tighten this bound to $70$.)
    \item For each $0\le j\le n-1$, we have $\abs{\sec_n(t^{-j}S t^j)}\le \abs S+10$.
    \item We have $\abs{\sec_n(U)} \le L$, for $L$ a constant depending on $S$, but crucially not on $n$.
\end{enumerate}
Combining these bounds with the previous observation, we would get
\[ \abs{ \sec_n\Big( N \cap\big((St)^*(St^{-1})^*\big)^m\Big)} \le (\abs{S}+10)\cdot L \cdot \left((\abs S+10)^{2n}\cdot 266\cdot L\right)^{m-1}\]
which is indeed (simply) exponential in $n$. Remains to prove estimates (a-b-c):

(a) Let $g\in\Gk$ and denote $h=t^{-n}gt^n$.
If $v=(0,\varepsilon)$, then $h_v\pmod{\Lc_n}$ is fully determined by the values $(\psi_{-2}(g),\psi_{-1}(g),\psi_0(g),\psi_1(g),\psi_2(g))$ (at most $8\cdot 4\cdot 2^3=256$ different sections).
Otherwise, $h_v$ is a subsection of some $\psi_*(g)\in \la a,d\ra$, i.e., belongs to $\{e,a,b,c,d,ad,da,ada,dad,adad\}$ ($10$ sections).

(b) Let $s\in S$ and denote $h=t^{-j}st^j$. If $v=(0,\varepsilon)$, then $h_v$ is fully determined by $s$ (so $\abs S$ choices). Otherwise, $h_v$ is a subsection of an element $\psi_*(g)\in \la a,d\ra$, i.e., belongs to $\{e,a,b,c,d,ad,da,ada,dad,adad\}$. We get $\abs S+10$ sections.
\begin{center}
\begin{minipage}{.48\linewidth}
    \centering
    \begin{tikzpicture}[scale=.8, xscale=-1]
        \clip (-1,-1) rectangle (8,9);
        \emptleft1{1.5}
        \node at (0,-.7) {\small$g$};
        \emptright1{1.5}
        \node at (2,-.75) {\footnotesize$\psi_{-1}(g)$};
        \emptright23
        \node at (3,.75) {\footnotesize$\psi_{-2}(g)$};
        \emptright3{4.5}
        \node at (4,2.25) {\footnotesize$\psi_{-3}(g)$};
        \emptright46
        \node at (5,3.75) {\footnotesize$\psi_{-4}(g)$};
        
        \emptright{5.5}{8.25}
        \node at (6.5,6) {\footnotesize$\psi_*(g)$};
        
        \draw (1,1.5) -- (4.33,6.5);
        \draw[dashed] (4.33,6.5) -- (5.13,7.75);
        \draw (5.13,7.75) -- (5.5,8.25);
        \draw[dashed] (5.5,8.25) -- (6,9);
        
        \draw[dotted] (-.5,8.25) -- (7.5,8.25);
        \draw[dotted] (-.5,0) -- (7.5,0);
        \node[fill=white, inner sep=2pt] at (6.7,0) {$\Lc_n$};
        
        \draw[fill=black] (5.5,8.25) circle (2pt);
        \draw[very thick, Green] (5.5,8.25) circle (5pt);
        \node[Green] at (5.5,8.7) {$v$};
    \end{tikzpicture}
    \captionof{figure}{Section of $t^{-n}gt^n$ below $v=(0,\varepsilon)$.}
\end{minipage}
\begin{minipage}{.48\linewidth}
\centering
        \begin{tikzpicture}[scale=.8, xscale=-1]
        \clip (-1,-1) rectangle (8,9);
        
        \draw (2,3) -- (4.33,-.5);
        \node at (2,1.1) {\large$s$};
        \emptright3{4.5}
        \node at (4,2.25) {\footnotesize$\psi_{-1}(s)$};
        \emptright46
        \node at (5,3.75) {\footnotesize$\psi_{-2}(s)$};
        
        \emptright{5.5}{8.25}
        \node at (6.5,6) {\footnotesize$\psi_*(s)$};
        
        \draw (-.33,-.5) -- (4.33,6.5);
        \draw[dashed] (4.33,6.5) -- (5.13,7.75);
        \draw (5.13,7.75) -- (5.5,8.25);
        \draw[dashed] (5.5,8.25) -- (6,9);
        
        \draw[dotted] (-.5,8.25) -- (7.5,8.25);
        \draw[dotted] (-.5,0) -- (7.5,0);
        \node[fill=white, inner sep=2pt] at (6.7,0) {$\Lc_n$};
        
        \draw[fill=black] (5.5,8.25) circle (2pt);
        \draw[very thick, Green] (5.5,8.25) circle (5pt);
        \node[Green] at (5.5,8.7) {$v$};
    \end{tikzpicture}
    \captionof{figure}{Section of $t^{-j}st^j$ below $v=(0,\varepsilon)$.}
\end{minipage}
\end{center}

(c) We'll need some more artillery this time. Let us start with a definition and some properties inspired by \cite{BOURABEE2010729}

\begin{defi}
Let $h\in N$. We define the \textbf{complexity level} $k(h)$ as the smallest integer $k$ such that all vertices $v\in\Lc_k$ have sections $h_v$ belonging to $\{e,a,b,c,d\}$ (with the convention $k(e)=k(b)=k(c)=k(d)=-\infty$).
\end{defi}
\begin{lemma}\label{lemma:complexity_levels}
For $g,h\in N$, we have
\begin{enumerate}[leftmargin=8mm, label={\normalfont(\roman*)}]
    \item $k(a)=0$
    \item $k(tgt^{-1})=k(g)-1$
    \item $k(gh)\le \max\{k(g),k(h)\}+1$.
\end{enumerate}
As a corollary, $k$ is finite on $N\setminus\{e,b,c,d\}$.
\end{lemma}
Let us go back to estimate (c). Let $K=\max_{s\in S}k(s)$. Each $h\in U$ can be written as
\[ h = ts_1t^{-1}\cdot t^2s_2t^{-2}\cdot\ldots \cdot t^Js_Jt^{-J}\cdot t^{J-1}s_{J+1}t^{1-J} \cdot \ldots \cdot ts_{2J-1}t^{-1} \]
with $m\ge 1$ and $s_1,s_2,\ldots,s_{2J-1}\in S$. Using Lemma \ref{lemma:complexity_levels} (ii-iii) repetitively gives $k(h)\le K+1$. It follows that, for any vertex $v\in\Lc_0$, the section $h_v$ is given by some automorphism of the $K+1$ first levels and then a choice from $\{e,a,b,c,d\}$ for the section at each vertex of the $(K+1)$-th level. In total, this gives a bound of
\[ L = 2^{2^{K+1}-1}\cdot 5^{2^{K+1}} \]
possibles sections.
\end{proof}

%\textbf{Question.} The complexity of the word problem in $\Gk$ has received some attention. However, the Dehn function has obviously not been investigated, as $\Gk$ is not finitely presented\footnote{Actually, we can define a Dehn function for any finitely generated group w.r.t.\ any given set of relations, even infinite. It just won't be a group invariant.}. What about the Dehn function for $\Lk$ ?

%\import{Source_files/}{sec8_questions.tex}

\printbibliography

\end{document}